\documentclass[final,12pt]{elsarticle}

\usepackage{appendix}
\usepackage{amssymb}
\usepackage{amsmath}
\usepackage{appendix}
\usepackage{graphicx}
\usepackage{caption}
\usepackage{subcaption}
\usepackage{color,soul}
\usepackage{url}
\usepackage{multirow}
\usepackage{float}
\usepackage{mathrsfs}
\usepackage[normalem]{ulem}

\graphicspath{{images/}}

\newcommand{\nospacetimes}{{\mkern-2mu\times\mkern-2mu}}

    \newtheorem{remark}{Remark}

\journal{Journal of Computational Physics}

\begin{document}

\begin{frontmatter}

\title{An isogeometric boundary element method for electromagnetic scattering with compatible B-spline discretizations}

\author[glasgow]{R. N. Simpson\corref{cor1}}
\ead{robert.simpson.2@glasgow.ac.uk}

\author[glasgow]{Z. Liu}

\author[lausanne,pavia]{R. V\'{a}zquez}

\author[ucb]{J.A. Evans}

\cortext[cor1]{Corresponding author}

\address[glasgow]{School of Engineering, University of Glasgow, Glasgow G12 8QQ, U.K.}

\address[lausanne]{Institute of Mathematics, \'Ecole Polytechnique F\'ed\'erale de Lausanne, Station 8, 1015 Lausanne, Switzerland}

\address[pavia]{Istituto di Matematica Applicata e Tecnologie Informatiche ``E.\ Magenes'' del CNR, via Ferrata 5, 27100, Pavia (Italy)}
\address[ucb]{Department of Aerospace Engineering Sciences, University of Colorado Boulder, Boulder, CO 80305, USA}

\begin{abstract}

We outline the construction of compatible B-splines on 3D surfaces that satisfy the continuity requirements for electromagnetic scattering analysis with the boundary element method (method of moments).  Our approach makes use of Non-Uniform Rational B-splines to represent model geometry and compatible B-splines to approximate the surface current, and adopts the isogeometric concept in which the basis for analysis is taken directly from CAD (geometry) data.  The approach allows for high-order approximations and crucially provides a direct link with CAD data structures that allows for efficient design workflows. After outlining the construction of div- and curl-conforming B-splines defined over 3D surfaces we describe their use with the electric and magnetic field integral equations using a Galerkin formulation.  We use B\'{e}zier extraction to accelerate the computation of NURBS and B-spline terms and employ $\mathscr{H}$-matrices to provide accelerated computations and memory reduction for the dense matrices that result from the boundary integral discretization.  The method is verified using the well known Mie scattering problem posed over a perfectly electrically conducting sphere and the classic NASA almond problem.  Finally, we demonstrate the ability of the approach to handle models with complex geometry directly from CAD without mesh generation. 

\end{abstract}


\end{frontmatter}

\section{Introduction}
\label{sec:introduction}

Research into unifying geometry and analysis for efficient design workflows has progressed rapidly in recent years driven  by the isogeometric analysis and computational geometry research communities.  Analysis based on geometry discretizations now covers a wide range of technologies including NURBS \cite{Hughes2005}, T-splines \cite{Bazilevs2010}, LR B-splines \cite{Johannessen2014}, PHT-splines \cite{wang2011adaptive} and subdivision surfaces \cite{cirak2000subdivision}.  A major research challenge at present is the automatic generation of volumetric discretizations from given geometric surface data and promising research includes the work of \cite{liu2015feature,wang2013trivariate} based on T-splines.  In contrast, analysis methods based on shell formulations or boundary integral methods are known to require only a surface discretization exhibiting key benefits for a common geometry and analysis model since no additional volumetric processing is required.  There has been much research into isogeometric shell formulations including \cite{cirak2000subdivision,benson2010isogeometric,Kiendl2009} and developments  into isogeometric boundary element methods based on NURBS \cite{Simpson2012, li2011isogeometric}, T-splines \cite{Scott2013, kostas2015ship} and subdivision surfaces \cite{bandara2015boundary}. 

A key application of the boundary element method is the analysis of electromagnetic scattering over complex geometries in which a perfectly electrically conducting (PEC) assumption can be made.  The method is often termed the method of moments within the electromagnetic research community but is synonymous with the Galerkin boundary element method.  It is well known that a straightforward application of nodal basis functions to the electric and magnetic field integral equations (EFIE, MFIE)  prevents numerical convergence and instead, discrete spaces that satisfy the relevant continuity requirements must be used.   The most commonly used discretization that satisifes the relevant continuity requirements are Raviart-Thomas \cite{raviart1977mixed} or RWG \cite{rao1982electromagnetic} basis functions that are mainly based on low order polynomials.  

In the context of isogeometric analysis progress has been made on the development of spline-based compatible discretizations \cite{buffa2011isogeometric,evans2013isogeometric, vazquez2010isogeometric, buffa2010isogeometric} in which a discrete de Rham sequence can be constructed providing a crucial step towards application of isogeometric analysis for fluid flow and electromagnetics applications.  This fundamental work opens up the opportunity for the development of an isogeometric boundary element method (isogeometric method of moments) for electromagnetic scattering which is the focus of the present study.  We note similar work in which subdivision surfaces are employed \cite{Li2016145}, but we believe that use of B-spline based algorithms provides greater refinement flexibility, provide a natural link with NURBS based systems that are ubiquitous in modern engineering design software, and offer higher convergence rates over equivalent subdivision schemes with extraordinary points.

We organise the paper as follows: first, we prescribe the Galerkin formulation of the relevant integral equations that govern electromagnetic scattering; we give an overview of NURBS surfaces and detail the construction of compatible B-splines; we then specify the fully discretized form of the integral equations for electromagnetic scattering with compatible B-splines; we cover implementation details of the method including fast evaluation of basis functions through B\'{e}zier extraction and the use of $\mathscr{H}$-matrices to approximate dense matrices; we verify the present method by performing electromagnetic scattering over a sphere in which a closed-form solution is provided by Mie scattering theory and finally, we demonstrate the ability of the present approach to perform electromagnetic scattering of PEC bodies with complex geometries taken directly from CAD software. It is assumed that time-harmonic fields are prescribed and, unless stated otherwise, it can be assumed that $\mathbf{x}\in \mathbb{R}^3$.

\section{Electric field integral equation: Galerkin formulation}
\label{sec:emagdiscretization}

\begin{figure}[h]
	\centering
	\includegraphics[width=0.7\textwidth]{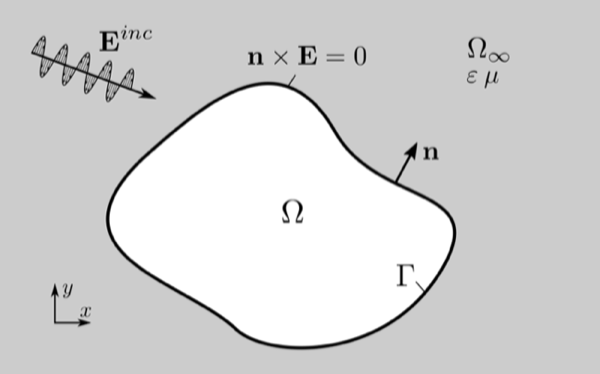}
	\caption{A PEC domain residing within an infinite domain impinged by an electromagnetic plane wave.}  
	\label{fig:d-pec-domain}
\end{figure}

We first assume a PEC domain $\Omega$ with connected boundary $\Gamma := \partial \Omega$ residing within an unbounded domain $\Omega_{\infty}$ with isotropic permeability and permittivity given by the scalar quantities $\varepsilon$ and $\mu$ respectively. We further assume a polarised time-harmonic electromagnetic plane wave of angular frequency $\omega$ is imposed on the PEC body with a wavenumber $k=\omega \sqrt{\varepsilon \mu}$.  Denoting $\mathbf{E}$ as the total electric field, in the presence of an electromagnetic wave a surface current $\mathbf{J}$ is induced and the following PEC condition holds on the surface of the scattered object
\begin{equation}
\mathbf{n} \times \mathbf{E} = 0 \label{eq:PECcondition1}
\end{equation}
where $\mathbf{n}$ represents the outward pointing normal vector.  We specify the incident wave as $\mathbf{E}^{i}(\mathbf{x}) = \mathbf{p} \, e^{-jk \mathbf{d}\cdot \mathbf{x}}$
where $j$ is the unit imaginary number, $\mathbf{p} = (p_x, p_y, p_z)$ is a polarization vector and $\mathbf{d} = (d_x, d_y, d_z), |\mathbf{d}| = 1$  is a propagation vector. The relationship between the total, incident and scattered electric fields is written as
\begin{equation}
  \mathbf{E} = \mathbf{E}^{i} + \mathbf{E}^s \label{eq:Etotal}\\
\end{equation}
where $\mathbf{E}^s$ represents the scattered electric field. The entire set-up is depicted in Figure~\ref{fig:d-pec-domain}.
%
%
%
%

Following the potential formulation of Maxwell's equations (see e.g. \cite{harrington1961time}), the scattered electric field can be expressed in terms of an  electric potential $\varphi$ and magnetic vector potential $\mathbf{A}$ (assuming time-harmonic fields) as
\begin{equation}	
  \label{eq:escattered_potentialformulation}
  \mathbf{E}^s = -j\omega \mathbf{A} - \nabla \varphi
\end{equation}
where the electric potential is given by 
\begin{equation}
  \label{eq:electric_potential}
  \varphi(\mathbf{x}) = \frac{1}{\varepsilon} \int_\Gamma \rho \, \frac{e^{-jkr}}{4 \pi r} \, \mathrm{d}\Gamma(\mathbf{y})
\end{equation}
with $r := |\mathbf{x} - \mathbf{y}|$ and the charge density $\rho$ expressed as
\begin{equation}
  \label{eq:charge_density}
  \rho = -\frac{1}{j\omega} \nabla \cdot \mathbf{J}
\end{equation}
with the magnetic potential related to the surface current through
\begin{equation}
  \label{eq:magnetic_potential}
  \mathbf{A}(\mathbf{x}) = \mu \int_\Gamma \mathbf{J}(\mathbf{y})\, \frac{e^{-jkr}}{4 \pi r} \, \mathrm{d}\Gamma(\mathbf{y}).
\end{equation}
We omit variable dependencies in future equations where they are implied by their context and adopt the notation $\Gamma_{y} \equiv \Gamma(\mathbf{y})$ and $\Gamma_{x} \equiv \Gamma(\mathbf{x})$. Substituting \eqref{eq:electric_potential} and \eqref{eq:magnetic_potential} into \eqref{eq:escattered_potentialformulation} and employing \eqref{eq:charge_density} with $k^{2}=\omega^{2} \varepsilon \mu$ and $j^{2}=-1$, the scattered electric field is expressed in terms of surface quantites as
\begin{equation}
  \label{eq:scattered_efield2}
  \mathbf{E}^s = - j\omega \mu \left( \int_{\Gamma_{y}} \mathbf{J} \frac{e^{-jkr}}{4 \pi r} \, \mathrm{d}\Gamma_{y} + \frac{1}{k^2} \nabla_{\Gamma_{x}} \int_{\Gamma_{y}} \nabla_{\Gamma_{y}} \cdot \mathbf{J}\, \frac{e^{-jkr}}{4 \pi r} \, \mathrm{d}\Gamma_{y} \right)
\end{equation}
where $\nabla_{\Gamma_{x}}$, $\nabla_{\Gamma_{y}}$ are surface gradient operators taken with respect to $\mathbf{x}$ and $\mathbf{y}$ respectively. Defining the linear operator 
\begin{equation}
  \label{eq:linear_operator}
  L^E[\boldsymbol{\tau}(\mathbf{x})] = \int_{\Gamma_{y}} \boldsymbol{\tau} \frac{e^{-jkr}}{4 \pi r} \, \mathrm{d}\Gamma_{y} + \frac{1}{k^2} \nabla_{\Gamma_{x}}  \int_{\Gamma_{y}} \nabla_{\Gamma_{y}}  \cdot \boldsymbol{\tau}\, \frac{e^{-jkr}}{4 \pi r} \, \mathrm{d}\Gamma_{y}
\end{equation}
along with the force term $\mathbf{f} = (j\omega \mu)^{-1} \mathbf{E}^i$, the Galerkin formulation of the EFIE reads as:

given $\mathbf{f}$, find $\mathbf{J} \in \mathcal{V}$ such that 
\begin{equation}
  \label{eq:galerkin_efie1}
   \langle \mathbf{w}, L^E[\mathbf{J}] \rangle  = \langle \mathbf{w}, \mathbf{f} \rangle \quad \forall \,\mathbf{w} \in \mathcal{V}
\end{equation}
where $\mathcal{V}$ is the trace space $H^{-\frac{1}{2}}(\mathrm{div}_\Gamma, \Gamma)$, and the $\langle \cdot, \cdot \rangle$ is the duality pairing between $\mathcal{V}$ and $H^{-\frac{1}{2}}(\mathrm{curl}_\Gamma, \Gamma)$. When the fields are smooth enough, the duality pairing reduces to $\langle \mathbf{u}, \mathbf{v} \rangle = \int_\Gamma \mathbf{u}\cdot \mathbf{v}\, \mathrm{d} \Gamma$. 

We define the finite dimensional subspace $\mathcal{V}_{h} \subset \mathcal{V}$ which allows the solution of \eqref{eq:galerkin_efie1} to be approximated as the solution of

given $\mathbf{f}$, find $\mathbf{J}_h \in \mathcal{V}_h$ such that
\begin{equation}
  \label{eq:galerkin_efie_h}
     \langle \mathbf{w}_h, L^E[\mathbf{J}_h] \rangle  = \langle \mathbf{w}_h, \mathbf{f} \rangle \quad \forall \, \mathbf{w}_h \in \mathcal{V}_h.
\end{equation}
Conventionally, $ \mathbf{w}_h$ and $\mathbf{J}_h$ are discretized through the Raviart-Thomas basis, but in our approach we make use of compatible B-splines that we now outline in detail.

\section{Discretization}
\label{sec:discretization}

\subsection{NURBS surfaces}
\label{sec:nurbs-surfaces}

Our implementation assumes a watertight NURBS surface parameterization that may be composed of multiple patches and we further assume that the connectivity of global basis functions between NURBS patches is known \emph{a priori}.  Dealing with the single patch case first, a NURBS surface parameterization is defined through a set of four-dimensional homogeneous control points $\{\mathbf{P}_{a}\}_{a=1}^{n_{p}}$, $\mathbf{P}_{a} = (x_{a}w_{a},y_{a}w_{a},z_{a}w_{a}, w_{a})$ (where $w_{a}$ represents a control point weight), a set of knot vectors $\{\Xi_{i}\}_{i=1}^{2}$ where $\Xi_{1} = \{0 = s_{1}, s_{2}, \ldots, s_{n+p+1} = 1\}$, $\Xi_{2} = \{0 = t_{1}, t_{2}, \ldots, t_{m+q+1} = 1\}$ and a degree vector $\mathbf{p}= (p,q)$. $n$ and $m$ denote the number of basis functions defined through the knot vectors $\Xi_1$ and $\Xi_2$ respectively with $n_{p} = n \times m$. We assume all knot vectors are open (i.e. for a given degree $p$ the knot vector contains $p+1$ equal knot values at its beginning and end). 

Defining the parametric domain $\widehat{\Gamma} = (0,1)^2 \subset \mathbb{R}^2$ and physical domain $\Gamma \subset \mathbb{R}^3$, a NURBS geometric mapping $\mathbf{F}: \widehat{\Gamma} \to \Gamma$ can be written in terms of parametric coordinates  $\mathbf{s} = (s,t) \in \widehat{\Gamma}$ as
\begin{equation}\label{eq:geometric-mapping-F}
\mathbf{F} = \sum_{a=1}^{n_{p}} R_{a}(\mathbf{s}) \mathbf{P}_{a}
\end{equation}
with the set of rational basis functions  $\{R_{a}\}_{a=1}^{n_{p}}$ defined as

%
%
\begin{equation}
R_{a}(\mathbf{s}) \equiv R_{a}(s,t) =  \frac{w_a B_a(s,t)}{\sum_{b=1}^{nm} w_b B_b(s,t)} \quad a=1,2,\ldots n_p
\end{equation}
where
\begin{equation}
B_a(s,t) = B_i^p(s)B_j^p(t), \nonumber
\end{equation}
with the set of univariate B-spline basis functions $\{B^{p}_i\}_{i=1}^n$ defined through the Cox-de-Boor algorithm (see e.g.~\cite{piegl2012nurbs}). The parametric basis function index $a$ is defined in terms of the univariate basis indices $i,j$ through
\begin{equation}\label{eq:parametric-basis-defn}
a = (j - 1)n  + i.
\end{equation}

Defining vectors of unique knot values in the $s$ and $t$ parametric directions as $\boldsymbol{\zeta}_1 = \{\zeta_1^1, \zeta_2^1, \ldots \zeta_{n_k}^1\}$ and $\boldsymbol{\zeta}_2 = \{\zeta_1^2, \zeta_2^2, \ldots \zeta_{m_k}^2\}$ respectively, the mesh in the parametric domain is given by
\begin{equation}
\mathcal{M}_h = \{ Q = (\zeta_i^1, \zeta_i^1+1) \times (\zeta_j^2, \zeta_j^2+1),\, 1\leq i \leq n_k - 1, \, 1 \leq j \leq m_k-1 \}
\end{equation}
with $n_e = \textrm{size}(\mathcal{M}_h)$ denoting the number of elements within the patch. Each element $Q$ within the patch contains $(p+1)\times(q+1)$ non-zero basis functions.


\subsection{Compatible B-spline approximation}

Given a set of univariate B-spline basis functions $\{B_i^{p}\}_{i=1}^n$, the space spanned by this basis is defined as
\begin{equation}\label{eq:univariate-bspline-span}
\widehat{S}^p := \textrm{span}\{B_i^{p}\}_{i=1}^n
\end{equation}
and in a similar manner, the tensor product B-spline space defined through the set of B-spline basis functions $B_{a} := B_i^{p} \otimes B_j^{q}$, $i=1,2,\ldots,n$, $j=1,2,\ldots m$ is defined as
\begin{equation}\label{eq:bspline-tensor-basis-space}
\widehat{S}^{p,q}:= \widehat{S}^p \otimes \widehat{S}^q = \textrm{span}\{B_a\}_{a=1}^{n_{b}}
\end{equation}
where the mapping defined by \eqref{eq:parametric-basis-defn} is employed and a hat symbol denotes that the quantity is defined over the parametric domain.  A div-conforming vector B-spline space is defined over the parametric domain as
\begin{equation}\label{eq:div-conforming-space-parametric}
\widehat{S}_1 := \widehat{S}^{p,q-1} \times \widehat{S}^{p-1,q}
\end{equation}
and likewise, a curl-conforming vector B-spline space is defined as
\begin{equation}\label{eq:curl-conforming-space-parametric}
\widehat{S}_2 := \widehat{S}^{p-1,q} \times \widehat{S}^{p,q-1}.
\end{equation}
The equivalent div-conforming and curl-conforming spaces defined in the physical domain are then constructed through appropriate Piola mappings as
\begin{equation}\label{eq:div-conforming-space-physical}
\mathcal{U}_h = \{  \mathbf{u} : \mathbf{u} \circ \mathbf{F} = \frac{1}{J} D\mathbf{F}\, \widehat{\mathbf{v}}, \, \widehat{\mathbf{v}} \in  \widehat{S}_1\}\\
\end{equation}
and
\begin{equation}\label{eq:curl-conforming-space-physical}
\mathcal{V}_h = \{ \mathbf{v} : \mathbf{v} \circ \mathbf{F} = \left(D\mathbf{F}^{+}\right)^{\mathrm{T}}\, \widehat{\mathbf{v}}, \, \widehat{\mathbf{v}} \in  \widehat{S}_2\}
\end{equation}
respectively, where $D\mathbf{F}$ is the Jacobian associated with the geometric mapping $\mathbf{F}$ which for 3D surfaces is given by the rectangular matrix
\begin{equation}\label{eq:rectangular-jacobian}
D\mathbf{F} = 
\begin{bmatrix}
\frac{\partial x}{\partial s} & \frac{\partial x}{\partial t} \\[0.3em]
\frac{\partial y}{\partial s} & \frac{\partial y}{\partial t} \\[0.3em]
\frac{\partial z}{\partial s} & \frac{\partial z}{\partial t}
\end{bmatrix},
\end{equation}
$D\mathbf{F}^{+}$ is the Monroe-Penrose pseudoinverse of the Jacobian given by
\begin{equation}\label{eq:jacobian-pseudoinverse}
D\mathbf{F}^{+} = \left( D\mathbf{F}^T D\mathbf{F} \right)^{-1} D\mathbf{F}^T,
\end{equation}
and $J$ is the surface element given by
\begin{equation}\label{eq:square-jacobian-Dvalue}
J = \sqrt{\left(\frac{\partial y}{\partial s}\frac{\partial z}{\partial t} - \frac{\partial z}{\partial s}\frac{\partial y}{\partial t}  \right)^2 +
	\left(\frac{\partial z}{\partial s}\frac{\partial x}{\partial t} - \frac{\partial x}{\partial s}\frac{\partial z}{\partial t}  \right)^2 +
	\left(\frac{\partial x}{\partial s}\frac{\partial y}{\partial t} - \frac{\partial y}{\partial s}\frac{\partial x}{\partial t}  \right)^2}.
\end{equation}
Further details of the derivation of \eqref{eq:div-conforming-space-physical} and \eqref{eq:curl-conforming-space-physical} can be found in \cite{evans2013isogeometric, buffa2010isogeometric} and the derivation of \eqref{eq:rectangular-jacobian}-\eqref{eq:square-jacobian-Dvalue} can be found in \cite[Sect.~5.4]{peterson2005mapped}.

\subsubsection{Basis functions}

Expressing vectors within the parametric domain as $\widehat{\mathbf{v}} = \widehat{v}_i \widehat{\mathbf{e}}_i$, $i=1,2$ and adopting the notation $\{B_a^{(p,q-1)}\}_{a=1}^{n_b^1}$, $\{B_a^{(p-1,q)}\}_{a=1}^{n_b^2}$ to represent the set of B-spline basis functions associated with the spaces $\widehat{S}^{p,q-1}$ and $\widehat{S}^{p-1,q}$ respectively, the set of div-conforming basis functions in the parametric domain $\hat \Gamma$ is defined as
\begin{equation}\label{eq:parametric-div-bspline-basis}
\widehat{\mathbf{N}}_a^{\textrm{div}}(s,t) = 
\begin{cases}
B_a^{(p,q-1)}(s,t)\,\widehat{\mathbf{e}}_1 \quad \quad \quad \,1 \leq a \leq n_b^1\\
B_{a - n_b^1}^{(p-1,q)}(s,t)\,\widehat{\mathbf{e}}_2 \quad  n_b^1+1 \leq a \leq n_b^1+n_b^2\\
\end{cases}
\end{equation}
which are transformed into a set of div-conforming basis functions on the surface $\Gamma$ using the Piola transformation defined in \eqref{eq:div-conforming-space-physical} as
\begin{equation}\label{eq:global-div-bspline-basis}
 \mathbf{N}_a^{\textrm{div}}(\mathbf{x}(s,t)) = \frac{1}{J} D\mathbf{F}\, \widehat{\mathbf{N}}_a^{\textrm{div}}(s,t) \quad 1 \leq a \leq n_b = n_b^1+n_b^2 
\end{equation}
where $\mathbf{F} \equiv \mathbf{F}(s,t)$ is implied. Curl-conforming basis functions are defined in analogous fashion.

\begin{figure}[htp]
	\centering
	\includegraphics[width=0.7\textwidth]{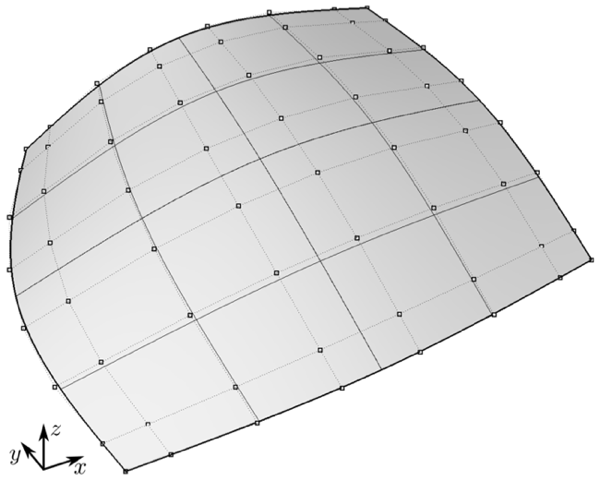}
	\caption{Bicubic NURBS patch defined by $n_{cp}=64$ control points, knot vectors $\Xi_1 = \Xi_2 =\{0,0,0,0,\frac{1}{4}, \frac{1}{2}, \frac{1}{2}, \frac{3}{4}, 1, 1, 1, 1\}$ and degrees $p=q=3$.  The degrees and knot vectors defined by the geometry are used directly to construct div-conforming B-splines.}
	\label{fig:nurbs-single-patch}
\end{figure}

\begin{figure}[htp]
	\centering
	\includegraphics[width=1.0\textwidth]{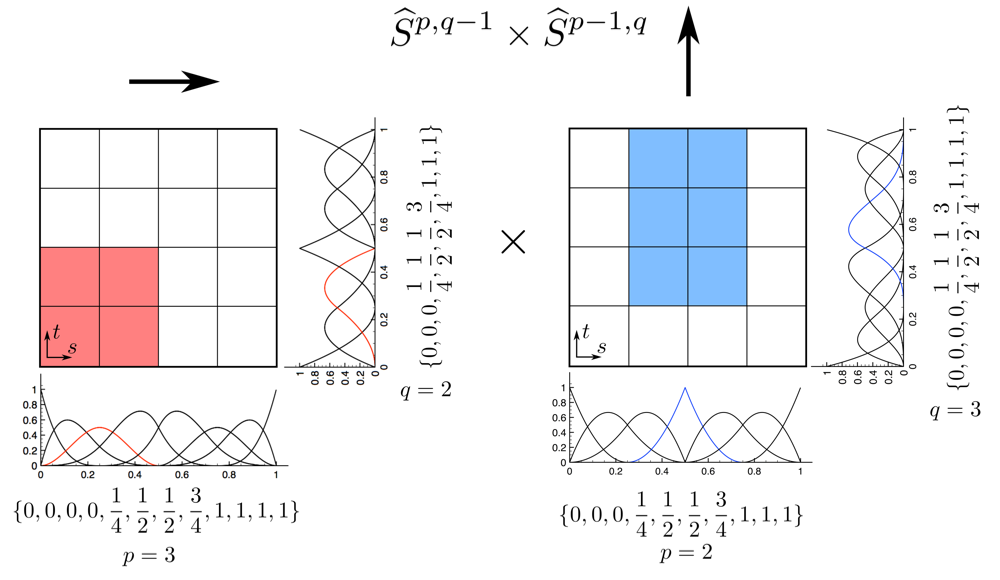}
	\caption{Construction of div-conforming basis functions defined over the parametric domain using the set of knot vectors and degrees defined by the geometry in Figure~\ref{fig:nurbs-single-patch}. The basis functions that define $\widehat{\mathbf{N}}_{19}^{\textrm{div}}(s,t)$ and the parametric interval that defines its span are highlighted in red.  Similarly for $\widehat{\mathbf{N}}_{88}^{\textrm{div}}(s,t)$ where all quantities are highlighted blue.}  
	\label{fig:div-basis-example-parametric}
\end{figure}

\begin{figure}
	\centering
	\begin{subfigure}[b]{0.6\textwidth}
		\includegraphics[width=\textwidth]{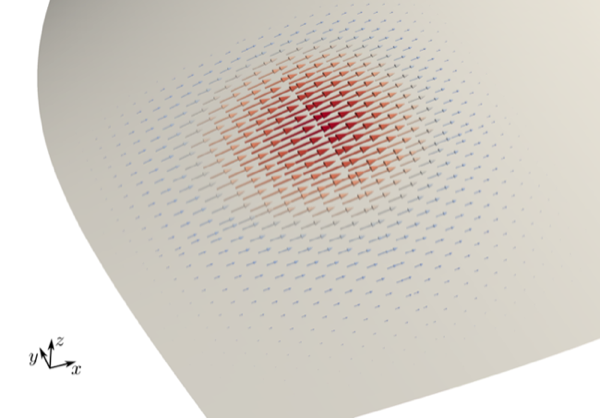}
		\caption{$\mathbf{N}_{19}^{\textrm{div}}(\mathbf{x})$}
		\label{subfig:div-conf-basis-uvector}
	\end{subfigure}
	
	\begin{subfigure}[b]{0.6\textwidth}
		\includegraphics[width=\textwidth]{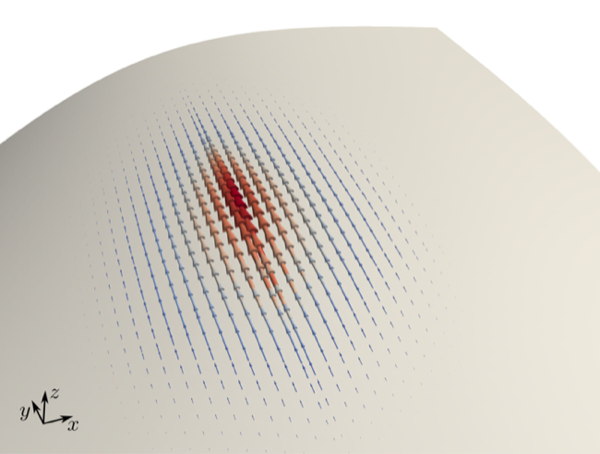}
		\caption{$\mathbf{N}_{88}^{\textrm{div}}(\mathbf{x})$}
		\label{subfig:div-conf-basis-vvector}
	\end{subfigure}
	\caption{Div-conforming B-splines defined over the surface given by the NURBS geometric mapping illustrated in Figure~\ref{fig:nurbs-single-patch}.  The basis functions correspond to those highlighted in Figure~\ref{fig:div-basis-example-parametric} where the Piola transform defined through \eqref{eq:div-conforming-space-physical} has been applied.}
	\label{fig:global-div-conf-basis}
\end{figure}

Global div- and curl-conforming approximations in physical space can then simply be expressed through
\begin{equation}
\mathbf{u}_h^{\textrm{div}}(\mathbf{x}) = \sum_{a=1}^{n_b} \mathbf{N}_a^{\textrm{div}}(\mathbf{x}) u_a
\end{equation}
and 
\begin{equation}
\mathbf{v}_h^{\textrm{curl}}(\mathbf{x}) = \sum_{a=1}^{n_b} \mathbf{N}_a^{\textrm{curl}}(\mathbf{x}) v_a
\end{equation}
respectively, where $u_a$ and $v_a$ are control coefficients. To illustrate the construction of compatible B-splines based on the NURBS parameterization shown in Figure~\ref{fig:nurbs-single-patch},  the bivariate B-splines generated from univariate B-splines are shown for two example basis functions in Figure~\ref{fig:div-basis-example-parametric}.  Further application of the Piola transformation as defined in \eqref{eq:global-div-bspline-basis} generates the div-conforming B-spline basis functions in physical space as shown in Figure~\ref{fig:global-div-conf-basis}.

\begin{remark}
For simplicity the construction of compatible B-splines is described using the same degree $(p, q)$ of the geometry. In practice it is possible to use a different degree for the B-splines discretization, as we will see in the numerical experiments.
\end{remark}

\subsection{Multipatch discretizations}
\label{sec:multipatch-discretization}

\begin{figure}
	\centering
	\begin{subfigure}[b]{0.6\textwidth}
		\includegraphics[width=\textwidth]{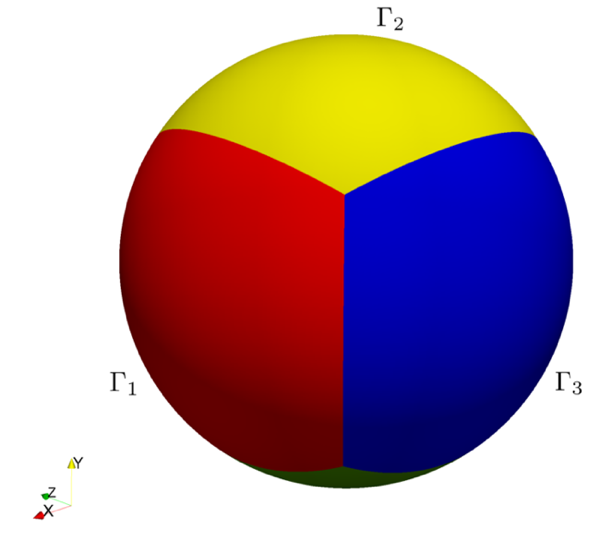}
		\caption{Physical domain $\Gamma$.}
		\label{subfig:nurbs-multipatch-physical}
	\end{subfigure}
	
	\begin{subfigure}[b]{0.6\textwidth}
		\includegraphics[width=\textwidth]{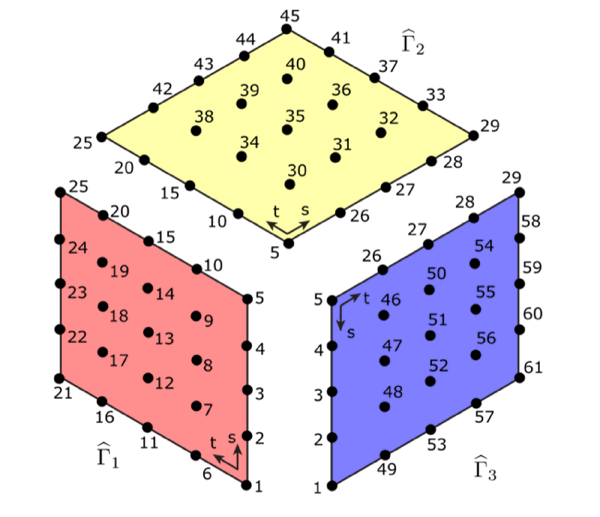}
		\caption{Parametric domains and geometry (nodal) connectivity.}
		\label{subfig:nurbs-multipatch-geom-connectivity}
	\end{subfigure}
	\caption{An example multipatch NURBS surface composed of patches of order $(4,4)$ with both physical and parametric domains illustrated.}
	\label{fig:nurbs-geometry-multipatch}
\end{figure}

Invariably, NURBS surfaces will consist of multiple patches whose union defines the physical domain through 
\begin{equation}
{\overline \Gamma} = \bigcup_{i=1}^{n_{d}} \overline{\Gamma_i}
\end{equation}
where $n_d$ is the number of parametric domains or patches and $\Gamma_i \cap  \Gamma_j = \emptyset$ for $i \neq j$. Each domain $\Gamma_i$ is constructed through a NURBS geometric mapping $\mathbf{F}_i: \widehat{\Gamma} \to \Gamma_i$ with parametric coordinates $\mathbf{s} \in \widehat{\Gamma}$ as
\begin{equation}\label{eq:geometric-mapping-F-multipatch}
\mathbf{F}_i = \sum_{a=1}^{n_{p}^i} R_{a}^i(\mathbf{s}) \mathbf{P}_{a}^i
\end{equation}
where the index $i$ indicates that the relevant quantity is restricted to patch $\Gamma_i$.  We require for two patches $\Gamma_i$ and $\Gamma_j$ with $i\neq j$ and which share a common edge the geometry mapping along the shared edge is the same. In addition, the knot vectors associated with each patch at the common edge must be the same, up to an affine transformation. Figure~\ref{subfig:nurbs-multipatch-physical} illustrates the geometry mappings of a multipatch NURBS geometry.

A global geometry connectivity array $C_g$ can be defined which maps a parametric basis function index $a$ and patch index $i$ to a global geometry basis index as
\begin{equation}\label{eq:geometry-conn}
A = C_g(i, a)  \quad  i=1,2\ldots n_d,\, a=1,2,\ldots n_p^i.
\end{equation}
The definition of the geometry connectivity array and the NURBS parameterisation given by \eqref{eq:geometric-mapping-F-multipatch} allows a multipatch NURBS parameterisation to be constructed such as that shown in Figure~\ref{subfig:nurbs-multipatch-geom-connectivity}.
%

As is well-known with vector bases, care must be taken when constructing global compatible basis functions since both the global basis function index and the orientation sign must be stored and we refer the reader to \cite{buffa2014isogeometric} where div- and curl-conforming B-spline approximations are constructed in a volumetric context.   We define the vector basis connectivity for a div-conforming basis through
\begin{align*}
A = C_{n}(i,a) \quad  i=1,2\ldots n_d,\, a=1,2,\ldots n_b^i
\end{align*}
where $n_b^i$ is the number of compatible B-spline basis functions in patch $i$. 
This allows a global multipatch compatible B-spline discretization to be written as
\begin{equation}\label{eq:conforming-mapping-multipatch}
\mathbf{u}^{\textrm{div}}_h(\mathbf{x})  =   \sum_{A=1}^{N_{b}} \mathbf{N}_A^{\textrm{div}}(\mathbf{x}) u_A
\end{equation}
where $N_b$ is the global number of basis functions, $\mathbf{N}_A^{\textrm{div}}|_{\Gamma_i} \equiv \mathbf{N}^{\textrm{div}}_{C_{n}(i,a)} \equiv \textrm{sgn}(i,a)\mathbf{N}^{\textrm{div}}_{i,a}$.

\begin{figure}
	\centering
	\begin{subfigure}[b]{0.45\textwidth}
		\includegraphics[width=\textwidth]{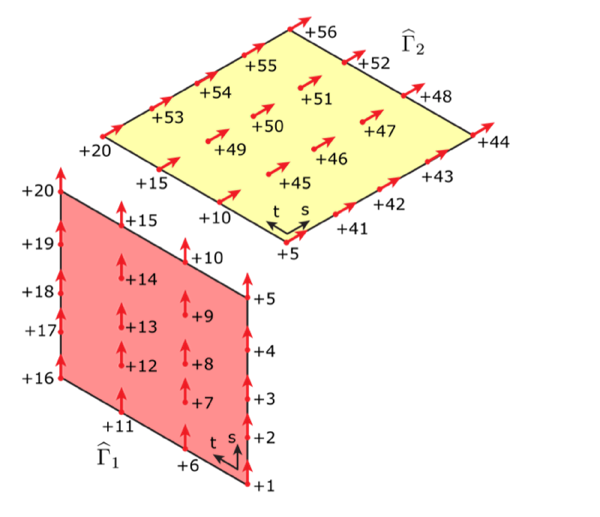}
		\caption{Domains $\widehat{\Gamma}_1$ and $\widehat{\Gamma}_2$.}
		\label{subfig:conforming-connectivitity-patches-12}
	\end{subfigure}
	\begin{subfigure}[b]{0.45\textwidth}
		\includegraphics[width=\textwidth]{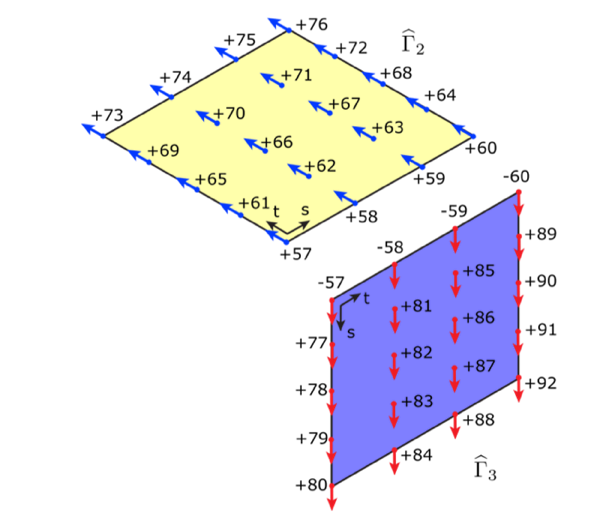}
		\caption{Domains $\widehat{\Gamma}_2$ and $\widehat{\Gamma}_3$.}
		\label{subfig:conforming-connectivitity-patches-23}
	\end{subfigure}

	\begin{subfigure}[b]{0.45\textwidth}
	\includegraphics[width=\textwidth]{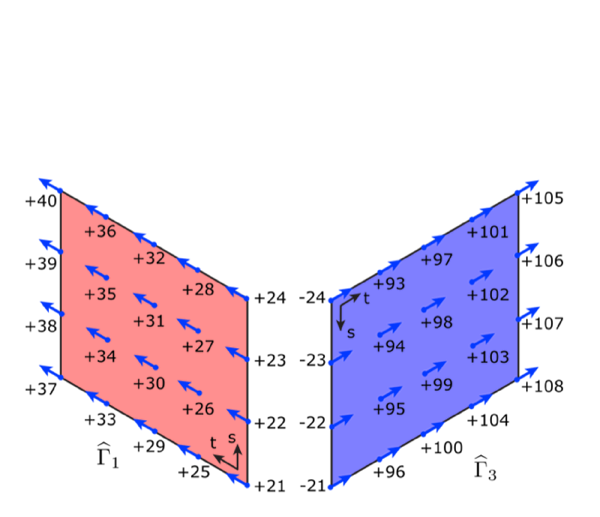}
	\caption{Domains $\widehat{\Gamma}_1$ and $\widehat{\Gamma}_3$.}
	\label{subfig:conforming-connectivitity-patches-13}
\end{subfigure}
	\caption{Example div-conforming vector basis connectivity associated with the NURBS multipatch geometry shown in Figure~\ref{fig:nurbs-geometry-multipatch} for a B-spline vector basis of order $(4,3)\nospacetimes(3,4)$. Red and blue arrows indicate a vector basis aligned in the $s$ and $t$ parametric directions respectively.}
	\label{fig:conforming-connectivity-multipatch}
\end{figure}

From an implementation standpoint the main consideration is how to handle basis functions along the edges of parametric domains which is best illustrated graphically.  Figure~\ref{fig:conforming-connectivity-multipatch} shows an example vector basis connectivity for div-conforming B-splines of order $(4,3)\nospacetimes(3,4)$ based on the geometry of Figure~\ref{fig:nurbs-geometry-multipatch}. Similar connectivities can be constructed for curl-conforming B-splines. 

\section{Discretised EFIE with compatible B-splines}

In the present work $\mathbf{w}_h$ and $\mathbf{J}_h$ in \eqref{eq:galerkin_efie_h} are defined through the the div-conforming B-spline discretization given by \eqref{eq:conforming-mapping-multipatch} and can be expressed as
\begin{align}
\mathbf{w}_h(\mathbf{x}) &= \sum_{A=1}^{N_b} \mathbf{N}^{\mathrm{div}}_A (\mathbf{x}) w_A \label{eq:weight_discretization}\\
 \mathbf{J}_h(\mathbf{x}) &=  \sum_{A=1}^{N_b} \mathbf{N}^{\mathrm{div}}_A (\mathbf{x}) j_A. \label{eq:current_discretization}
\end{align}
Substituting \eqref{eq:weight_discretization} and \eqref{eq:current_discretization} into \eqref{eq:galerkin_efie_h} and applying the divergence theorem to transfer a derivative onto $\mathbf{w}_h$, a system of equations is formed as
\begin{equation}
  \label{eq:system_equations_efie}
  \mathbf{Z}_{AB} \mathbf{J}_B = \mathbf{f}_A
\end{equation}
where 
\begin{equation}
    \begin{split}
      \label{eq:discretised_efie}
      \mathbf{Z}_{AB} = \int_{\Gamma_x} \mathbf{N}^{\mathrm{div}}_A \cdot &\left( \int_{\Gamma_y} \mathbf{N}^{\mathrm{div}}_B \frac{e^{-jkr}}{4 \pi r} \, \mathrm{d}\Gamma_{y} \right) \mathrm{d}\Gamma_{x} \\
        &- \frac{1}{k^2} \int_{\Gamma_x} \nabla_{\Gamma_{x}} \cdot \mathbf{N}^{\mathrm{div}}_A \left( \int_{\Gamma_y} \nabla_{\Gamma_{y}} \cdot \mathbf{N}^{\mathrm{div}}_B\, \frac{e^{-jkr}}{4 \pi r} \, \mathrm{d}\Gamma_{y} \right)\mathrm{d}\Gamma_{x}
    \end{split}
\end{equation}
\begin{equation}
    \mathbf{f}_A = \frac{1}{j\omega\mu}\int_{\Gamma_x} \mathbf{N}^{\mathrm{div}}_A \cdot \mathbf{E}^i\, \mathrm{d}\Gamma_{x}
\end{equation}
and $\mathbf{J}_B$ represents a vector of unknown surface current density coefficients. A similar procedure can be carried out for the magnetic field integral equation as detailed in Appendix~\ref{app:mfie}.

\subsection{Radar Cross Section}
\label{subsec:rcs}

The radar cross section $\sigma$ which quantifies how detectable an object is to a radar signal in a given direction is computed as
\begin{equation}
\label{eq:rcs}
\sigma = \lim_{R \to \infty} 4 \pi R^2 \frac{|\mathbf{E}^s|^{2}}{|\mathbf{E}^i|^{2}}
\end{equation}
where $R$ is the distance between the radar signal and the target object and furthermore, it can be assumed in the present work that $|\mathbf{E}^i|$ = 1. As detailed in \cite{gibson2008method, balanis2012advanced} if the source and field points are located far apart then  $R \approx |\mathbf{x}|$ and the scattered electric field at a source (observation) point can be expressed as
\begin{equation}\label{eq:efield-far-scattering}
\mathbf{E}^s(\mathbf{x})  = -\frac{j\omega \mu}{4 \pi} \frac{e^{-jk|\mathbf{x}|}}{|\mathbf{x}|} \int_{\Gamma_y} \mathbf{J}(\mathbf{y}) e^{j k \mathbf{d}\cdot \mathbf{y}} \, \mathrm{d}\Gamma_y
\end{equation}
allowing the RCS to be computed as
\begin{equation}
\label{eq:rcs2}
\sigma = 4 \pi |\mathbf{x}|^2|\mathbf{E}^s|^{2}
\end{equation}
or, in terms of the RCS in decibels per square metre
\begin{equation}
\label{eq:rcs-dbsm}
\sigma_{dBsm} = 10 \log_{10} \sigma.
\end{equation}

\section{Implementation}
\label{sec:implementation}

\begin{figure}[htp]
	\centering
	\includegraphics[width=0.7\textwidth]{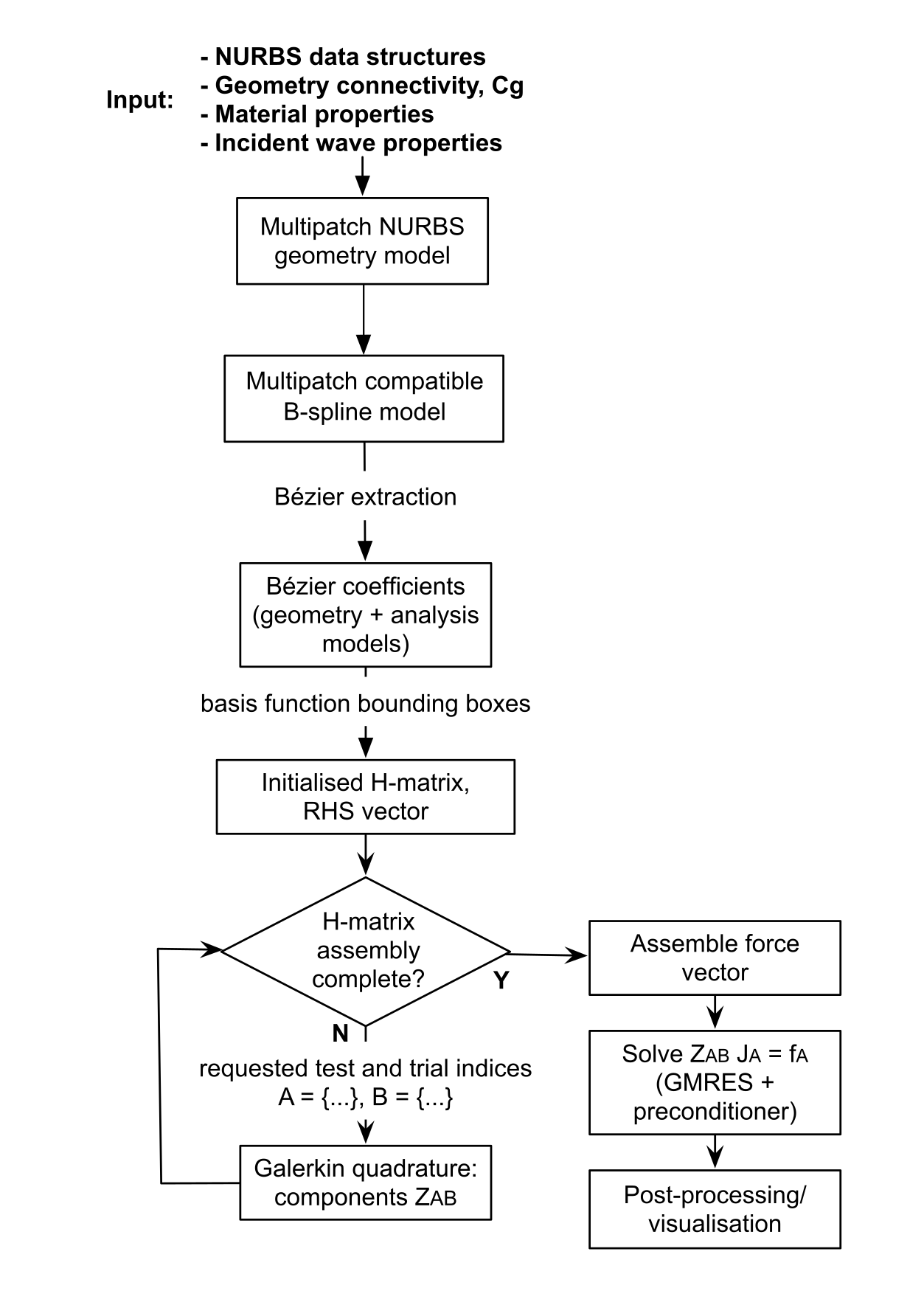}
	\caption{An outline of the algorithm for performing electromagnetic scattering with compatible B-splines using the boundary element method with $\mathscr{H}$-matrix acceleration and B\'{e}zier extraction.}
	\label{fig:algorithm-workflow}
\end{figure}

Figure~\ref{fig:algorithm-workflow} details the main steps in the implementation of the present method. A multipatch compatible B-spline discretization is constructed directly from the NURBS surface parameterization. The inherent link between the geometry and analysis models allows for straightforward computation of compatible basis functions with the relevant Piola transforms. We utilise B\'{e}zier extraction \cite{borden2011isogeometric} to accelerate computations whereby high order B-spline and NURBS basis functions are computed through precomputed B\'{e}zier extraction coefficients and inexpensive Bernstein polynomials.  

As is well-known with Galerkin boundary element methods, careful consideration must be given to the computation of the matrix components $\mathbf{Z}_{AB}$ given by \eqref{eq:discretised_efie} when the element domains $\Gamma_{x}$ and $\Gamma_{y}$ are either coincident, edge adjacent, vertex adjacent or lie close to one another.  We use the robust quadrature algorithms proposed by Sauter and Schwab \cite{sauter2010boundary} that deal with each of these cases.

\begin{figure}[htp]
	\centering
	\includegraphics[width=0.7\textwidth]{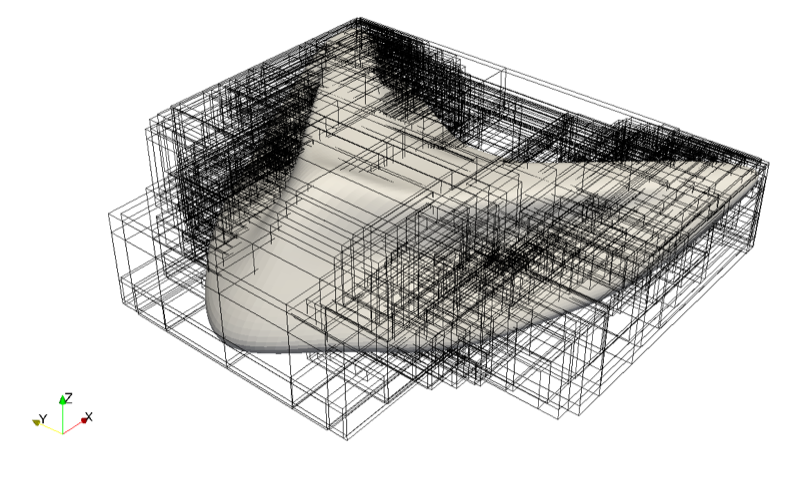}
	\caption{Example geometry with the corresponding set of bounding boxes defined by the support of each basis function used for low rank $\mathscr{H}$-matrix approximations.}
	\label{fig:example-bounding-boxes}
\end{figure}

To overcome the debilitating nature of large dense matrix $\mathbf{Z}$, we approximate this matrix using $\mathscr{H}$-matrices whereby a low-rank approximation is constructed through appropriate geometrical cluster trees that separate terms into admissible and non-admissible terms (i.e. far-field and near-field terms respectively).  We do not wish to delve into the technical details of $\mathscr{H}$-matrices and instead guide the reader to relevant literature (see e.g. \cite{bebendorf2008hierarchical, hackbusch2004hierarchical}).  However, we remark that $\mathscr{H}$-matrices are found to be particularly amenable for implementation into an existing BEM library and we make use of the library HLibPro \cite{kriemann2008hlibpro} which provides high-performance $\mathscr{H}$-matrix libraries that scale optimally over multicore hardware and are primarily based on the Adaptive Cross Approximation algorithm \cite{bebendorf2000approximation}.  The library requires as an input the set of bounding boxes defined by the support of each basis function (see Figure~\ref{fig:example-bounding-boxes}) and the basis function index associated with each box.  Once an $\mathscr{H}$-matrix approximation is formed for a particular wavenumber, the matrix can be written and read freely from file which allows for highly efficient radar cross section computations.  We note that this approach is valid for low to medium wavenumbers with special techniques required for high wavenumbers (e.g. \cite{bebendorf2015wideband}).

\section{Numerical results}
\label{sec:numericalresults}

To verify the present approach and to demonstrate the capability of the method of performing electromagnetic scattering directly from CAD models using an isogeometric approach we present numerical results for a series of electromagnetic scattering problems with PEC conditions. 

\subsection{PEC sphere}
\label{subsec:miescattering_sphere}

\begin{figure}[htp]
  \centering
\includegraphics[width=0.9\textwidth]{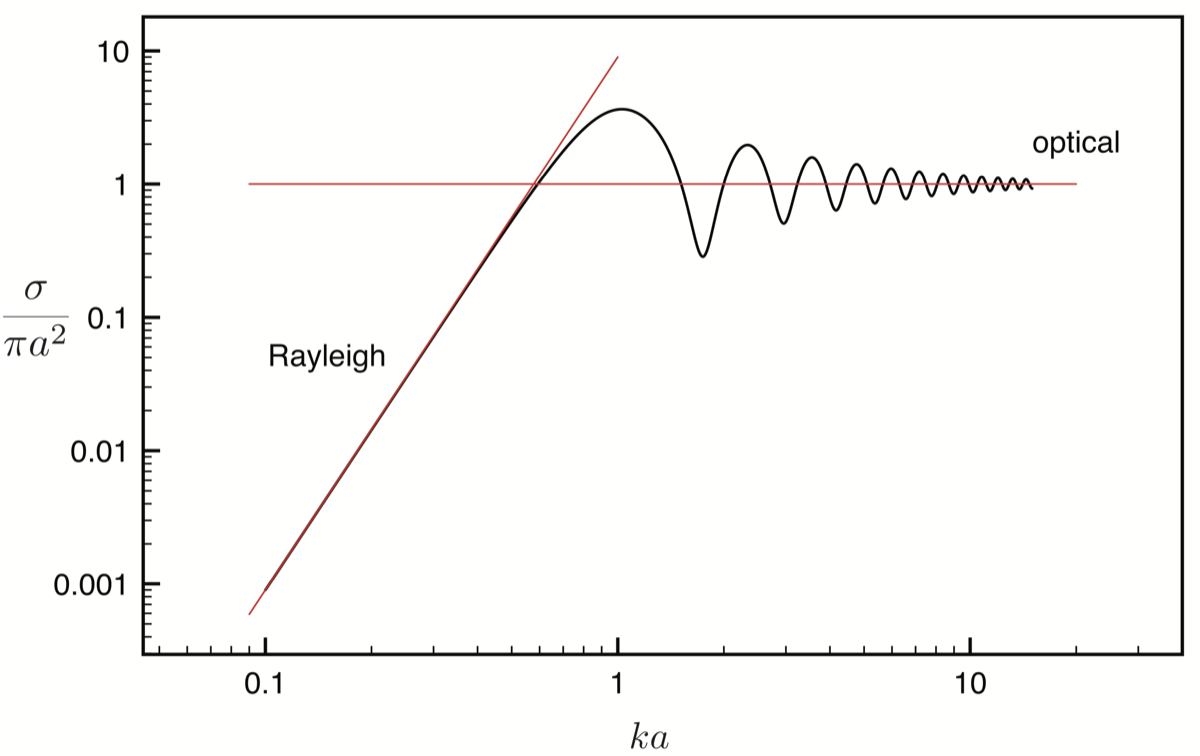}
  \caption{The monostatic radar cross section of a PEC sphere as a function of normalised wavenumber, commonly referred to as the Mie solution.}  
  \label{fig:mie-scattering-analytical}
\end{figure}

\begin{figure}[htp]
  \centering
\includegraphics[width=0.5\textwidth]{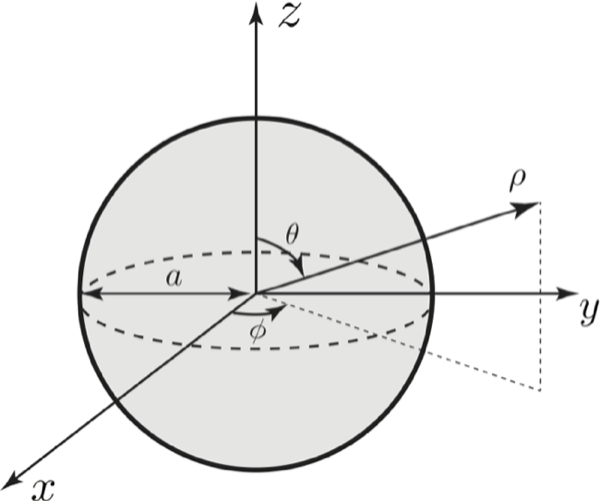}
  \caption{Mie scattering problem: spherical coordinate system.}  
  \label{fig:spherical-coord-system}
\end{figure}

The first problem we consider is that of electromagnetic plane wave impinging on a PEC sphere of radius $a=1$ which has a well-known solution given by the Mie series 
(see e.g. \cite{harrington1961time}).  The incident wave is polarised in the x-direction by specifying $\mathbf{p} = (1,0,0)$ and is chosen to propagate in the positive z-direction with $\mathbf{d} = (0, 0, 1)$.  The solution for the surface current given in spherical coordinates $(\rho, \theta, \phi)$ (see Figure~\ref{fig:spherical-coord-system}) is expressed as
\begin{align*}
  \label{eq:mie_scattering}
  J_{\rho} =& 0\\
  J_\theta =& \frac{j}{\eta} p_x \frac{\cos \phi}{ka} \sum_{n=1}^\infty a_n \left( \frac{\sin \theta P_n^{1'}(\cos \theta)}{\hat{H}_n^{(2)'}(ka)}  + \frac{j  P_n^{1}(\cos \theta)}{\sin \theta \hat{H}_n^{(2)}(ka)} \right)\\
  J_\phi =& \frac{j}{\eta} p_x \frac{\sin \phi}{ka} \sum_{n=1}^\infty a_n \left(  \frac{P_n^{1}(\cos \theta)}{\sin \theta \hat{H}_n^{(2)'}(ka)} - \frac{\sin \theta P_n^{1'}(\cos \theta)}{j\hat{H}_n^{(2)}(ka)} \right)
\end{align*}
with
\begin{equation}
a_n = \frac{j^{-n} (2n + 1)}{n(n+1)}
\end{equation}
where $\eta = \sqrt{\mu / \varepsilon}$, the terms $P_n^{1}$ and $P_n^{1'}$ correspond to the set of order 1 associated Legendre polynomials and derivatives respectively and  
\begin{align}
  \label{eq:spherical_hat_hankel1}
  \hat{H}_n^{(2)} &= k h^{(2)}_n\\
 \hat{H}_n^{(2)'} &= \left( nh^{(2)}_n - k h^{(2)}_{n+1} \right) + h^{(2)}_n
\end{align}
with $h^{(2)}_n$ denoting the spherical Hankel function of the second kind.  The radar cross section for this problem given in terms of increasing normalised wavenumber is illustrated in Figure~\ref{fig:mie-scattering-analytical} where the two asymptotic limits associated with Rayleigh and optical scattering are labelled.

Using the present approach, the sphere geometry is discretised using bi-quartic NURBS patches arranged in a cube topology with no degenerate points, as in Figure~\ref{subfig:nurbs-multipatch-physical}. Control point coordinates, weights and knot vectors for this NURBS parameterization can be found in \cite{cobb1988tiling}. We construct div-conforming B-splines using the knots inherited by the NURBS parameterization with degrees $(4,3){\nospacetimes}(3,4)$, $(3,2){\nospacetimes}(2,3)$, $(2,1){\nospacetimes}(1,2)$ and $(1,0){\nospacetimes}(0,1)$ and apply successive h-refinement (knot insertion) to generate a set of meshes h0 (base mesh), h1, h2 etc.  Table~\ref{table:mie-sphere-mesh-details} provides further details of each discretization. It should be noted that compatible B-splines of degree $(1,0){\nospacetimes}(0,1)$ are directly equivalent to low order Raviart-Thomas or RWG basis functions on quadrilateral meshes. The bi-quartic NURBS representation of the geometry is used for all analyses and thus geometric error is eliminated for all discretizations considered.

\begin{table}[]
\centering
\caption{Details of div-conforming B-spline discretizations used in the Mie scattering study.}
\label{table:mie-sphere-mesh-details}
\begin{tabular}{ccccc}
\hline
\multirow{2}{*}{\textbf{\begin{tabular}[c]{@{}c@{}}mesh \\ (\# elements)\end{tabular}}} & \multicolumn{4}{c}{\textbf{degrees of freedom}}                    \\ 
                                                                                        & $(1,0){\nospacetimes}(0,1)$& $(2,1){\nospacetimes}(1,2)$ & $(3,2){\nospacetimes}(2,3)$ & $(4,3){\nospacetimes}(3,4)$\\ \hline
h0 (6)	& 	12	&	48	&	108	&	192	\\
h1 (24)	&	48	&	108	&	192	&	300	\\
h2 (96)	&	192	&	300	&	432	&	588	\\
h3 (384)	& 	768	&	972	&	1,200	&	1,452	\\
h4 (1536)	&	3,072	&	3,468	&	3,888	&	4,332	\\ \hline
\end{tabular}
\end{table}

\begin{figure}[htp]
  \centering
\includegraphics[width=0.9\textwidth]{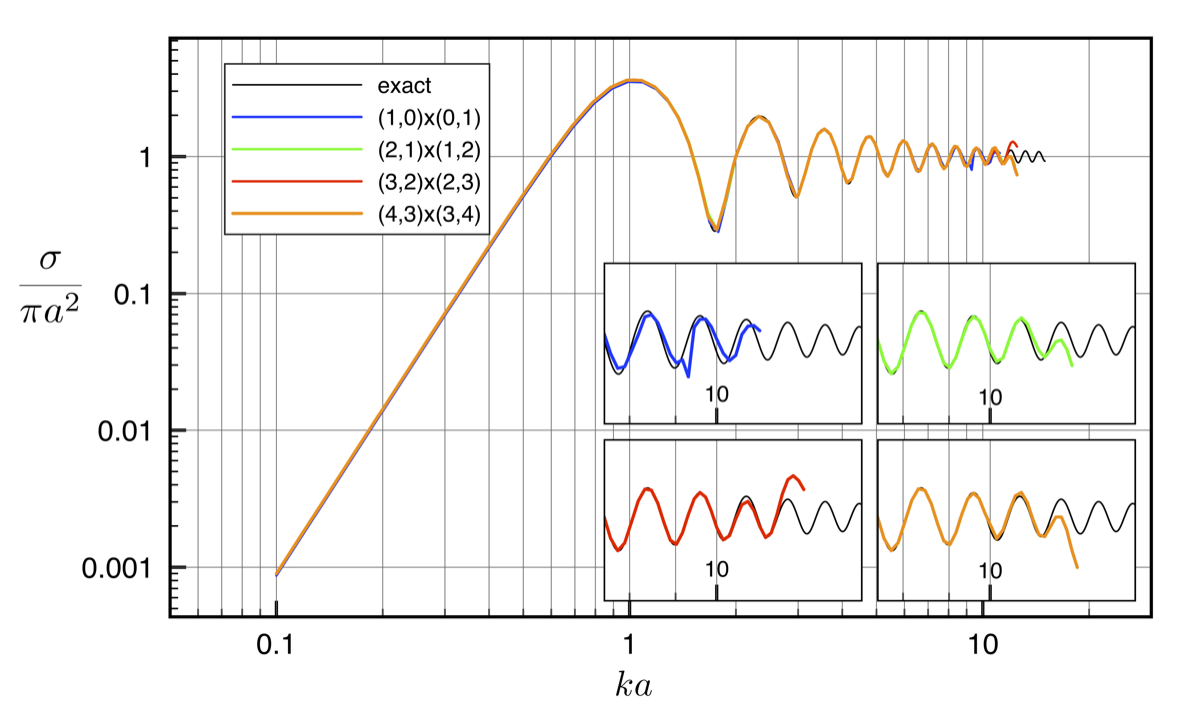}
  \caption{Normalised RCS values for a PEC sphere computed for increasing wavenumber with div-conforming B-splines of varying degree, mesh h3}  
  \label{fig:rcs-sphere-h3}
\end{figure}

After solving for surface current, equations \eqref{eq:efield-far-scattering} and \eqref{eq:rcs2} were used to determine radar cross section values with the results for mesh h3 shown in Figure~\ref{fig:rcs-sphere-h3} for each B-spline degree. The superior RCS accuracy obtained through higher order B-spline discretizations is demonstrated and this is also apparent in RCS values obtained with meshes h0, h1 and h2 as presented in Appendix~\ref{app:sec:pec-additional-results}.  As expected, finer meshes are capable of handling higher wavenumbers.  

\begin{figure}[htp]
	\centering
	\begin{subfigure}[b]{0.45\textwidth}
		\includegraphics[width=1.0\textwidth]{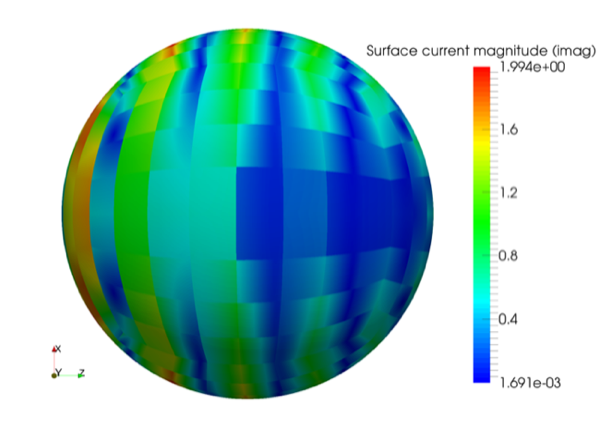}
		\caption{$|\mathbf{J}_i|$}
		\label{subfig:miesphere-surface-current-constant}
	\end{subfigure}
	\begin{subfigure}[b]{0.45\textwidth}
		\includegraphics[width=1.0\textwidth]{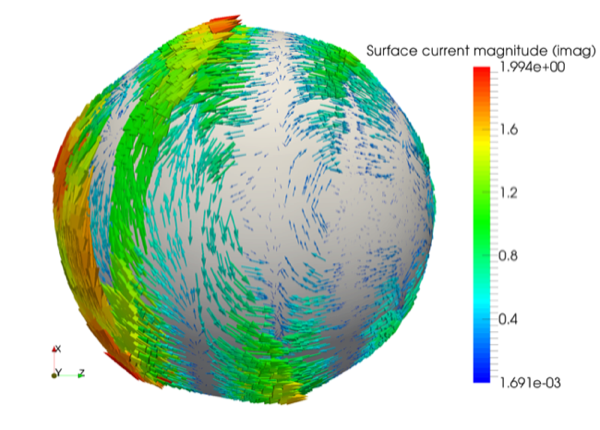}
		\caption{$\mathbf{J}_i$}
		\label{subfig:miesphere-surface-vector-constant}  
	\end{subfigure}
	
	\vspace{1ex}
	\begin{subfigure}[b]{0.8\textwidth}
	\includegraphics[width=1.0\textwidth]{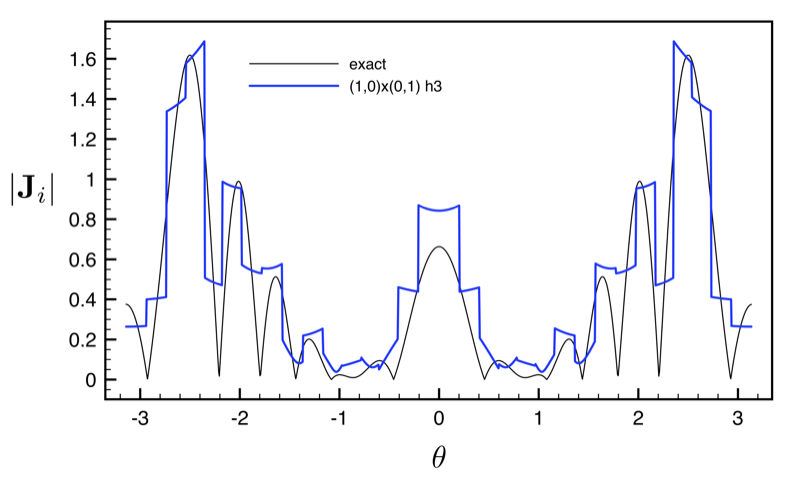}
	\caption{$|\mathbf{J}_i|$ sampled over $y\operatorname{-}z$ plane.}
	\label{subfig:miesphere-constant-sample-slice}  
	\end{subfigure}
	\caption{Sphere scattering problem, $k=8$: surface current quantities (imaginary component) obtained with div-conforming B-splines of degree $(1,0){\nospacetimes}(0,1)$ and three levels of h-refinement (h3).}
	\label{fig:miesphere-constant-h3}
\end{figure}

\begin{figure}[htp]
	\centering
	\begin{subfigure}[b]{0.45\textwidth}
		\includegraphics[width=1.0\textwidth]{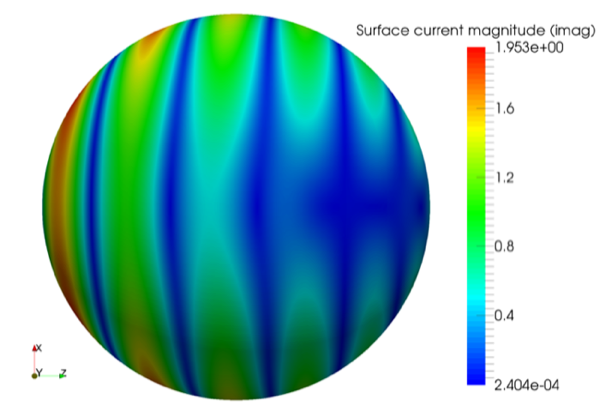}
		\caption{$|\mathbf{J}_i|$}
		\label{subfig:miesphere-surface-current-cubic}
	\end{subfigure}
	\begin{subfigure}[b]{0.45\textwidth}
		\includegraphics[width=1.0\textwidth]{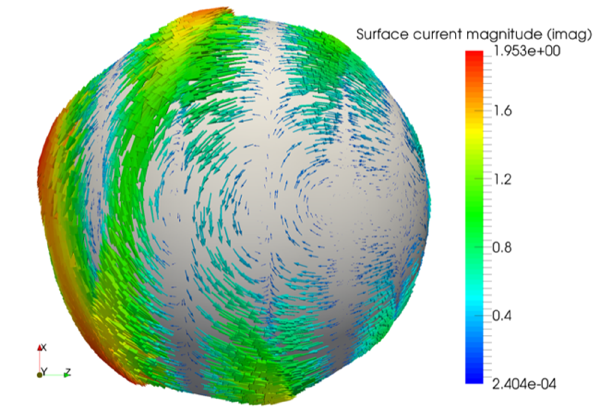}
		\caption{$\mathbf{J}_i$}
		\label{subfig:miesphere-surface-vector-cubic}  
	\end{subfigure}

	\vspace{1ex}
	\begin{subfigure}[b]{0.8\textwidth}
	\includegraphics[width=1.0\textwidth]{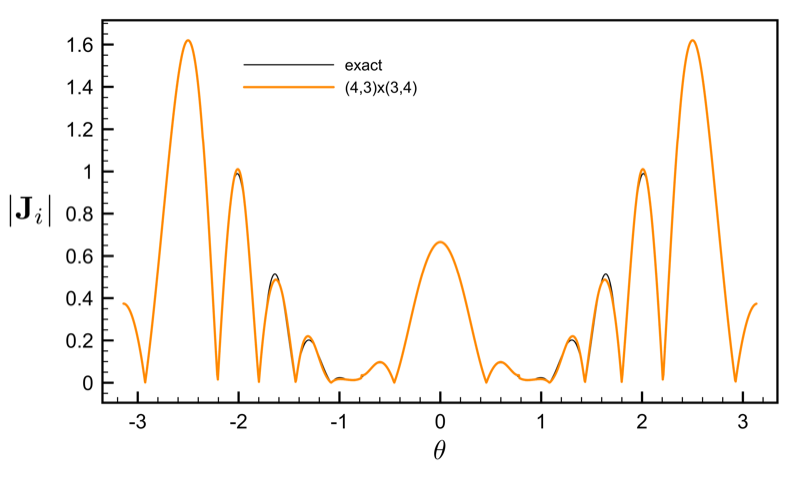}
	\caption{$|\mathbf{J}_i|$ sampled over $y\operatorname{-}z$ plane.}
	\label{subfig:miesphere-quadratic-sample-slice}  
	\end{subfigure}
	\caption{Sphere scattering problem $k=8$: surface current quantities (imaginary component) obtained with div-conforming B-splines of degree $(4,3){\nospacetimes}(3,4)$ and three levels of h-refinement (h3).}
	\label{fig:miesphere-cubic-h3}
\end{figure}

Plots of surface currents and magnitudes for $k=8$, h3 are shown for each B-spline degree in Figures~\ref{fig:miesphere-constant-h3} through to \ref{fig:miesphere-cubic-h3}  where the higher accuracy and smoothness offered through higher B-spline degrees is visible.  Recalling that the $(1,0){\nospacetimes}(0,1)$ discretization is equivalent to the commonly used Raviat-Thomas elements, it is clear that higher order compatible B-spline discretizations offer substantial accuracy improvements over such elements.

\begin{figure}[htp]
	\centering
	\begin{subfigure}[b]{0.45\textwidth}
		\includegraphics[width=1.0\textwidth]{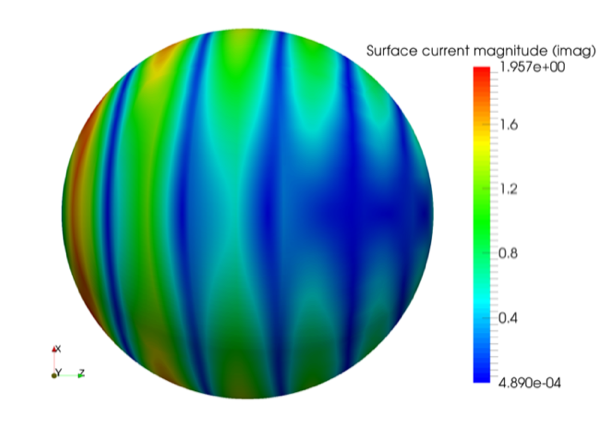}
		\caption{$|\mathbf{J}_i|$}
		\label{subfig:miesphere-surface-current-linear}
	\end{subfigure}
	\begin{subfigure}[b]{0.45\textwidth}
		\includegraphics[width=1.0\textwidth]{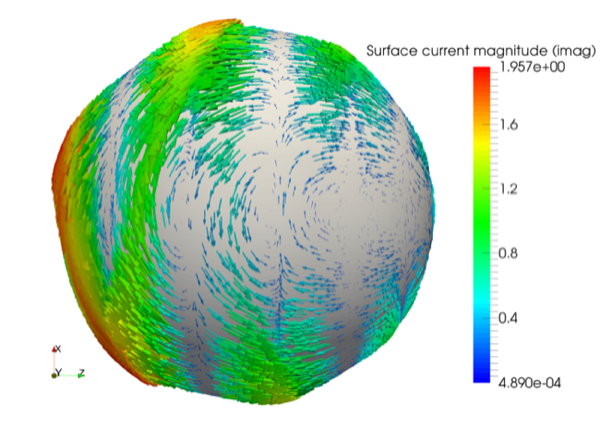}
		\caption{$\mathbf{J}_i$}
		\label{subfig:miesphere-surface-vector-linear}  
	\end{subfigure}

	\vspace{1ex}
	\begin{subfigure}[b]{0.8\textwidth}
	\includegraphics[width=1.0\textwidth]{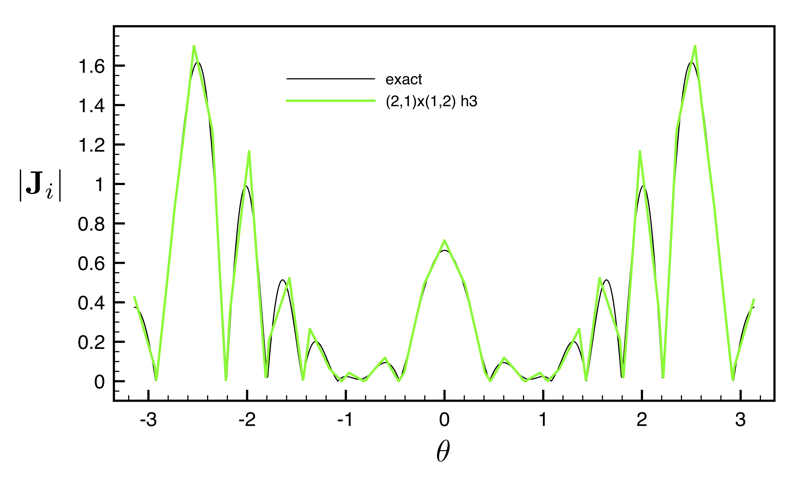}
	\caption{$|\mathbf{J}_i|$ sampled over $y\operatorname{-}z$ plane.}
	\label{subfig:miesphere-linear-sample-slice}  
	\end{subfigure}
	\caption{Sphere scattering problem $k=8$: surface current quantities (imaginary component) obtained with div-conforming B-splines of degree $(2,1){\nospacetimes}(1,2)$ and three levels of h-refinement (h3).}
	\label{fig:miesphere-linear-h3}
\end{figure}

\begin{figure}[htp]
	\centering
	\begin{subfigure}[b]{0.45\textwidth}
		\includegraphics[width=1.0\textwidth]{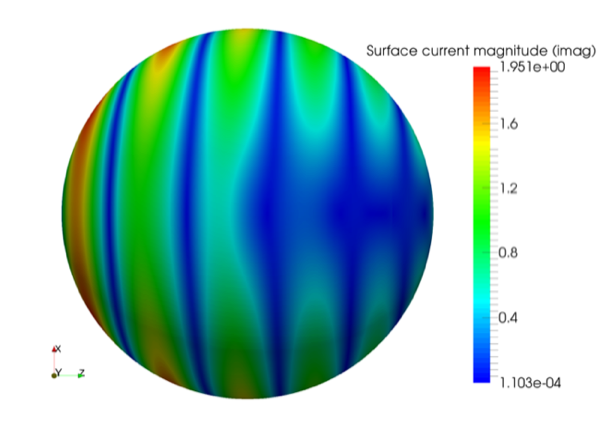}
		\caption{$|\mathbf{J}_i|$}
		\label{subfig:miesphere-surface-current-quadratic}
	\end{subfigure}
	\begin{subfigure}[b]{0.45\textwidth}
		\includegraphics[width=1.0\textwidth]{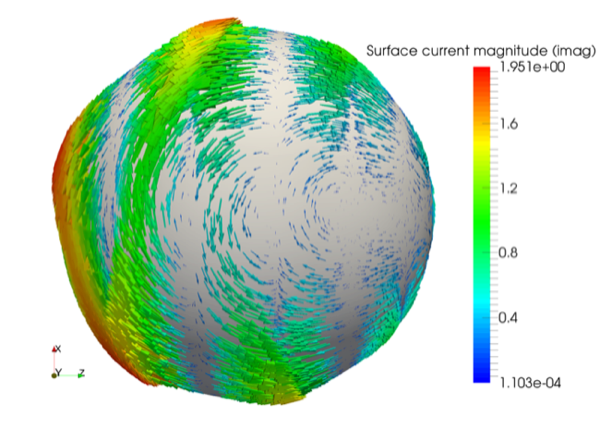}
		\caption{$\mathbf{J}_i$}
		\label{subfig:miesphere-surface-vector-quadratic}  
	\end{subfigure}

	\vspace{1ex}
	\begin{subfigure}[b]{0.8\textwidth}
	\includegraphics[width=1.0\textwidth]{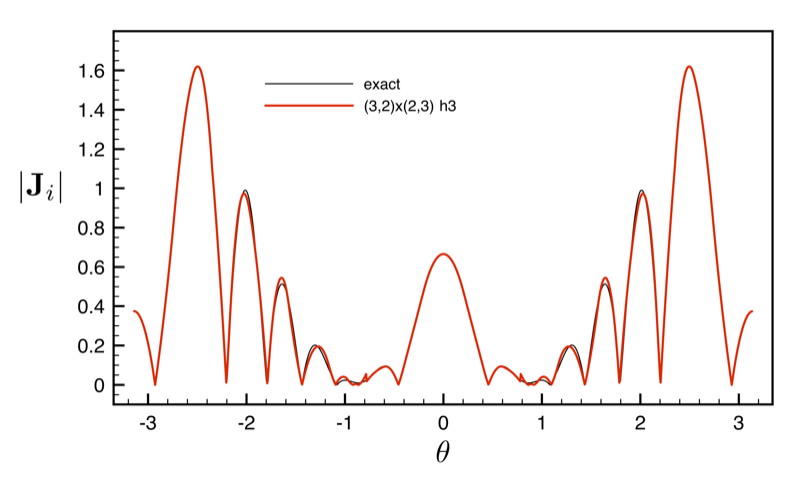}
	\caption{$|\mathbf{J}_i|$ sampled over $y\operatorname{-}z$ plane.}
	\label{subfig:miesphere-quadratic-sample-slice}  
	\end{subfigure}
	\caption{Sphere scattering problem $k=8$: surface current quantities (imaginary component) obtained with div-conforming B-splines of degree $(3,2){\nospacetimes}(2,3)$ and three levels of h-refinement (h3).}
	\label{fig:miesphere-quadratic-h3}
\end{figure}

\begin{figure}[htp]
	\centering
	\includegraphics[width=0.9\textwidth]{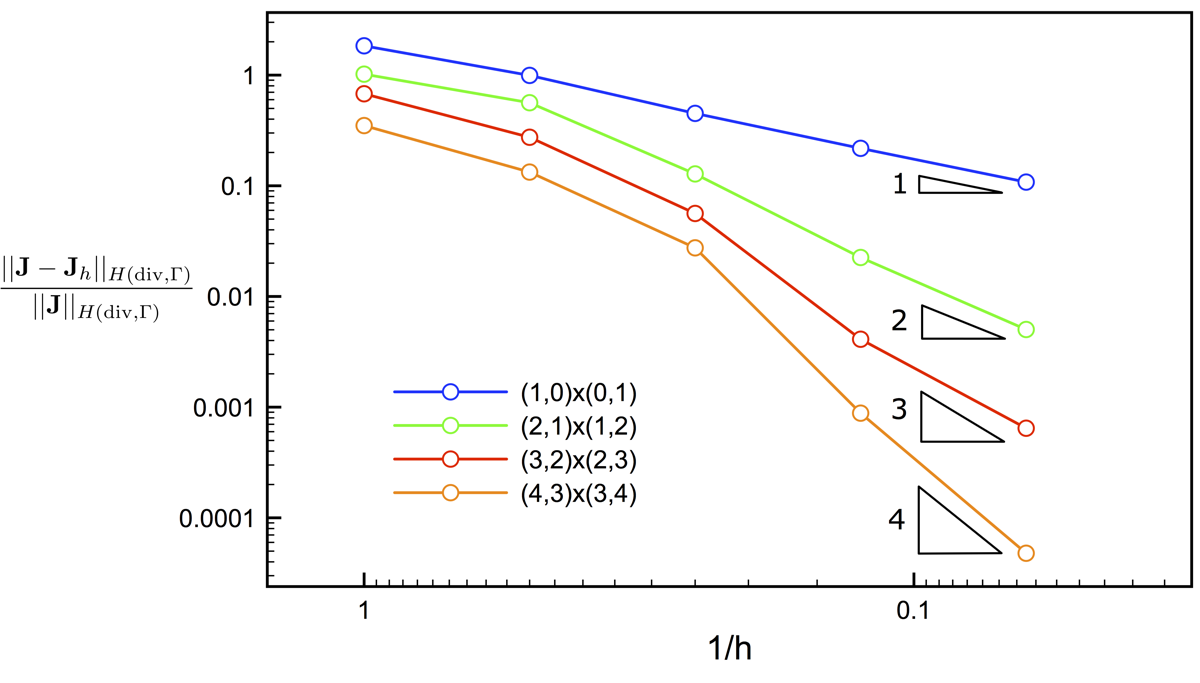}
	\caption{Mie scattering convergence study with $k = 3$: relative error norms for B-spline discretizations and theoretical convergence rates.}  
	\label{fig:mie-convergence-study}
\end{figure}

Additionally, to establish that correct convergence rates are obtained using our approach we compute relative errors using the norm defined by
\begin{equation}
  \label{eq:hdivnorm}
  || \mathbf{v} ||_{H(\text{div}, \Gamma)} = || \mathbf{v} ||_{L_{2}} +  || {\rm div}_{\Gamma}\,\mathbf{v} ||_{L_{2}},
\end{equation}
where we remark that the $L^2$ norm of the surface divergence is well defined for this particular example. A convergence rate of $p+1$ is expected for a given B-spline degree with minimum degree $p$.  We specify a wavenumber of $k=3$ and evaluate relative errors through the norm of \eqref{eq:hdivnorm} for each B-spline degree for meshes h0 to h4. Relative errors for this study are plotted in Figure~\ref{fig:mie-convergence-study} where theoretical convergence rates are demonstrated.
%
%
%

\subsection{NASA almond}
\label{sec:nasa_almond}

\begin{figure}[htp]
  \centering
\includegraphics[width=0.7\textwidth]{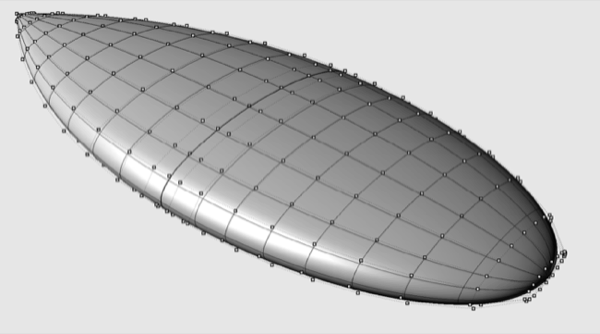}
  \caption{The NASA almond geometry represented by four bicubic NURBS patches with two degenerate points.}  
  \label{fig:nasa-almond-geometry}
\end{figure}

A common benchmark problem used to verify electromagnetic scattering numerical methods is the NASA almond problem as detailed in \cite{woo1993programmer}.   The geometry of the surface is defined through parametric expressions which are detailed in Appendix~\ref{app:sec:nasa-geometry}. In the present study these expressions were used as inputs to the Math Rhino plugin developed by Rhino3DE \cite{rhinoplugin} generating a NURBS representation of the almond geometry with four bicubic NURBS patches as shown in Figure~\ref{fig:nasa-almond-geometry}.  In addition, the software library Open CASCADE \cite{occt} was used to extract the necessary geometry data structures required to construct compatible B-spline discretizations defined over the almond surface. Div-conforming B-splines of orders $(3,2)\nospacetimes(2,3)$, $(2,1)\nospacetimes(1,2)$ and $(1,0)\nospacetimes(0,1)$ were generated with uniform h-refinement (knot insertion) applied to the initial discretization shown in Figure~\ref{fig:nasa-almond-geometry} to generate successively refined discretizations. Again, we use the notation h0, h1, h2 to indicate a mesh with no-refinement (base mesh), 1 level of h-refinement etc. and the abbreviations HH and VV to denote horizontally polarised and vertically polarised incident waves respectively. Table~\ref{table:nasa-almond-mesh-details} provides further details of each B-spline discretization. For the computation of the integrals we increase the number of quadrature points in the vicinity of the two degenerate points, to increase the accuracy.

\begin{table}[]
\centering
\caption{Details of compatible B-spline discretizations used for the NASA almond study.}
\label{table:nasa-almond-mesh-details}
\begin{tabular}{cccc}
\hline
\multirow{2}{*}{\textbf{\begin{tabular}[c]{@{}c@{}}mesh \\ (\# elements)\end{tabular}}} & \multicolumn{3}{c}{\textbf{degrees of freedom}}                    \\ 
                                                                                        & $(1,0){\nospacetimes}(0,1)$& $(2,1){\nospacetimes}(1,2)$ & $(3,2){\nospacetimes}(2,3)$\\ \hline
h0 (288)                                                                                & 558                  & 700                  & 858                  \\
h1 (1152)                                                                               & 2,268                & 2,546                & 2,840                \\
h2 (4608)                                                                               & 9,144                & 9,694                & 10,260               \\ \hline
\end{tabular}
\end{table}

To verify our implementation we compute the RCS given by \eqref{eq:rcs-dbsm} at frequencies of 1.19GHz, 3GHz and 7GHz for both horizontally and vertically polarised incident waves.  We use numerical RCS reference values from \cite{ganesh2006spectrally} for the 1.19GHz case, \cite{ganesh2006spectrally,antilla1994radiation} for the 3GHz case and \cite{feko} for the 7GHz case.  In addition, we utilise experimental results for the 1.19GHz  case as shown in \cite{woo1993programmer}. Both \cite{ganesh2006spectrally} and \cite{feko} are based on a boundary element (method of moments) approach with the work of \cite{antilla1994radiation} adopting a coupled finite element/boundary element formulation.

\begin{figure}[htp]
	\centering
	\includegraphics[width=1.0\textwidth]{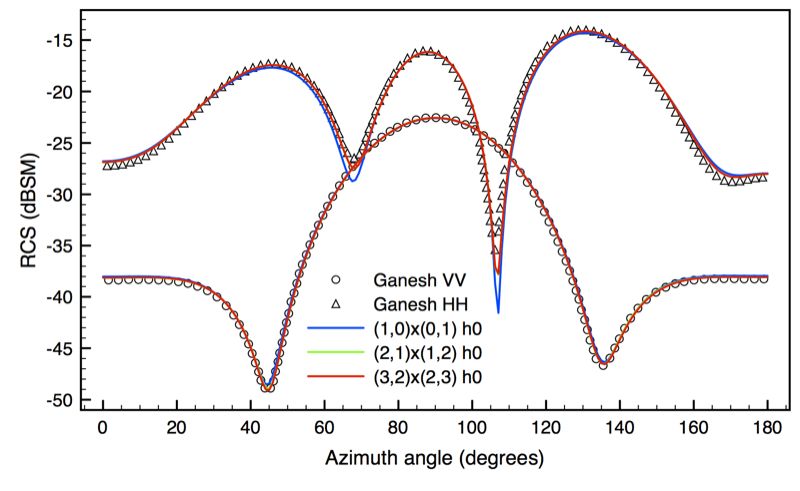}
	\caption{Radar cross section profile for NASA almond geometry: 1.19GHz, horizontal and vertical polarization.  Reference data obtained from \cite{ganesh2006spectrally}.}
	\label{fig:almond-1-19ghz-rcs}  
\end{figure}

\begin{figure}[htp]
	\centering
	\includegraphics[width=1.0\textwidth]{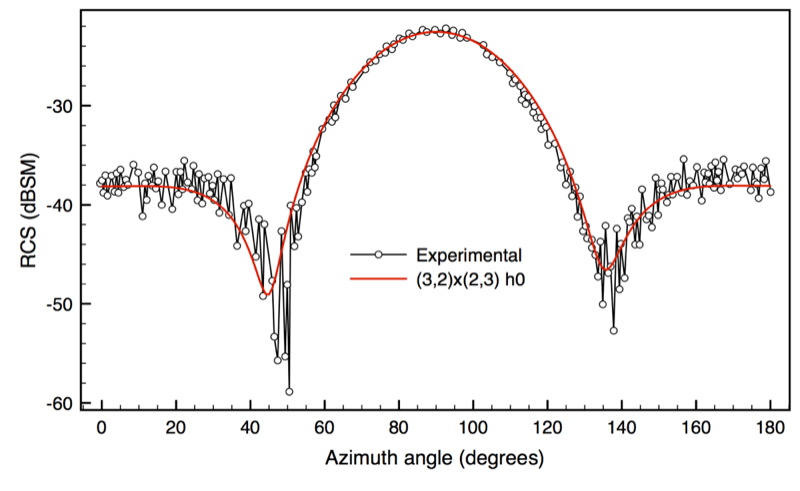}
	\caption{Comparison of experimental and numerical radar cross section profile for NASA almond geometry: 1.19GHz vertical polarization.  Experimental reference data obtained from \cite{woo1993programmer}.}
	\label{fig:almond-1-19ghz-rcs-experimental}  
\end{figure}

\begin{figure}[htp]
  \centering
\includegraphics[width=1.0\textwidth]{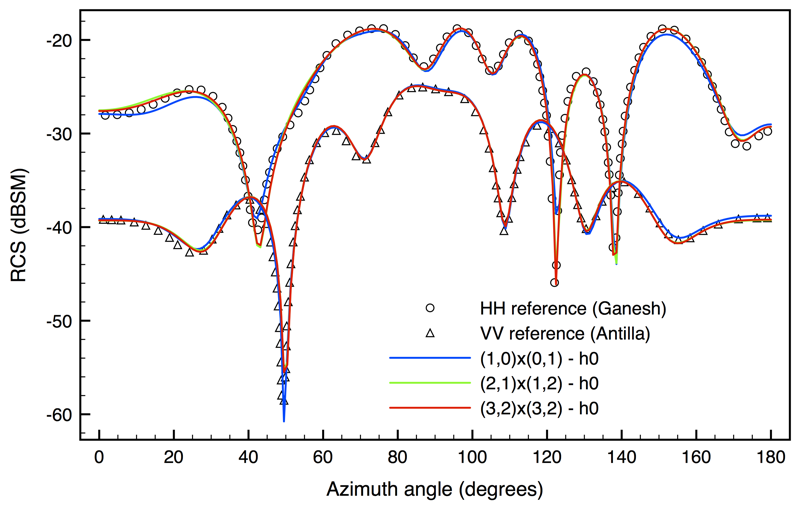}
  \caption{Radar cross section profile for NASA almond geometry: 3GHz, horizontal and vertical polarization.  Reference data obtained from \cite{ganesh2006spectrally,antilla1994radiation}.}
  \label{fig:almond-3ghz-rcs}  
\end{figure}

\begin{figure}[htp]
	\centering
	\begin{subfigure}[b]{0.6\textwidth}
		\includegraphics[width=1.0\textwidth]{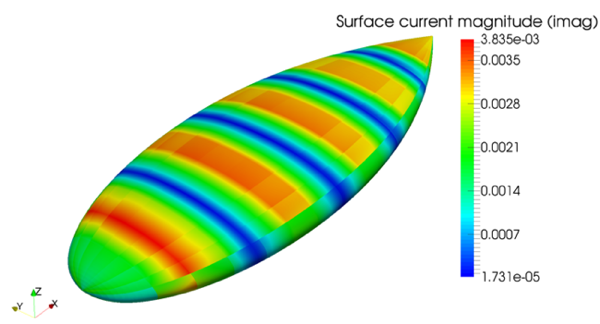}
		\caption{Order $(1,0)\nospacetimes(0,1)$.}
		\label{subfig:almond-3ghz-surface-current-loworder}
	\end{subfigure}
	
	\begin{subfigure}[b]{0.6\textwidth}
		\includegraphics[width=1.0\textwidth]{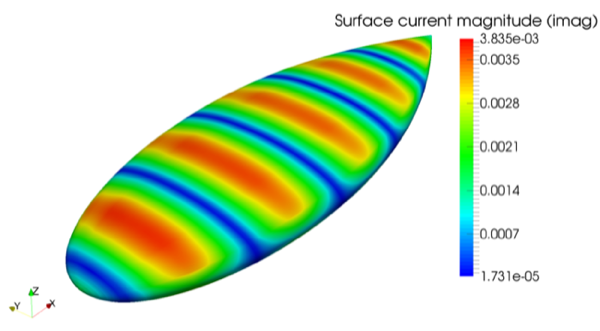}
		\caption{Order $(2,1)\nospacetimes(1,2)$.}
		\label{subfig:almond-3ghz-surface-current-midorder}  
	\end{subfigure}

	\begin{subfigure}[b]{0.6\textwidth}
		\includegraphics[width=1.0\textwidth]{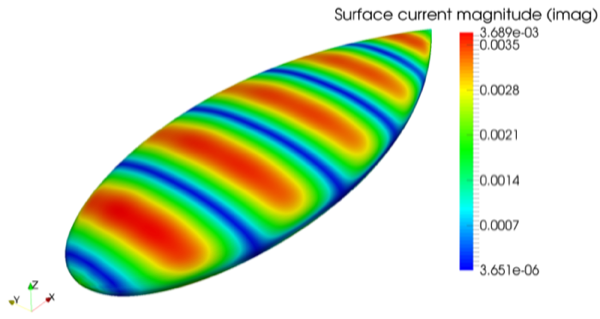}
		\caption{Order $(3,2)\nospacetimes(2,3)$.}
		\label{subfig:almond-3ghz-surface-current-highorder}  
	\end{subfigure}
	\caption{Magnitude of imaginary component of surface current over the NASA almond geometry: vertically polarised planewave of 3GHz travelling in the positive $x-$direction, mesh h0.}
	\label{fig:almond-3ghz-surface-current}
\end{figure}

%

\begin{figure}[htp]
  \centering
\includegraphics[width=1.0\textwidth]{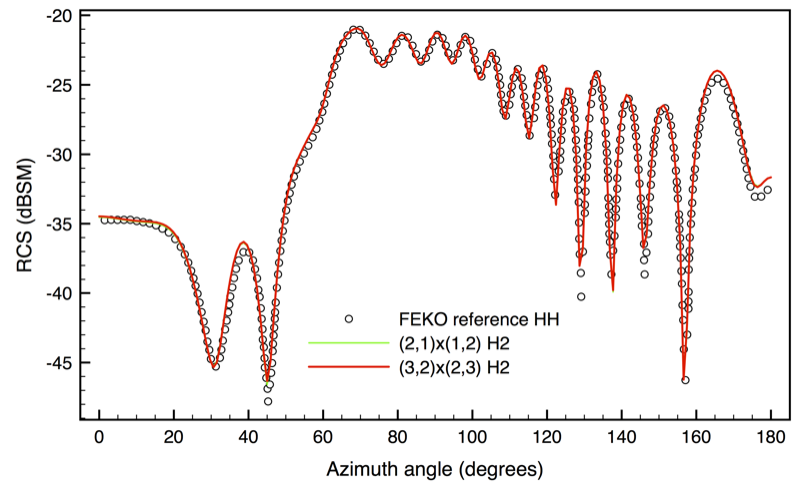}
  \caption{Radar cross section profile for NASA almond geometry: 7GHz, horizontal polarization.  Reference data obtained from \cite{feko}.}
  \label{fig:almond-7ghz-hh-rcs}  
\end{figure}

\begin{figure}[htp]
  \centering
\includegraphics[width=1.0\textwidth]{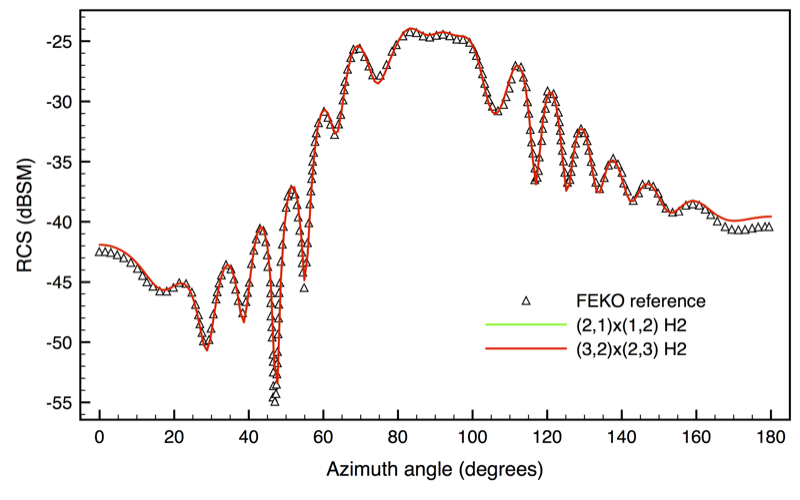}
  \caption{Radar cross section profile for NASA almond geometry: 7GHz, vertical polarization.  Reference data obtained from \cite{feko}.}
  \label{fig:almond-7ghz-vv-rcs}  
\end{figure}

Figure~\ref{fig:almond-1-19ghz-rcs} illustrates RCS plots for the 1.19GHz case for each B-spline order with mesh h0. Good agreement with the numerical reference solution is visible for each order. In addition, Figure~\ref{fig:almond-1-19ghz-rcs-experimental} demonstrates good agreement with experimental data for this frequency.  In a similar manner, numerical RCS values for the 3GHz case are shown in Figure~\ref{fig:almond-3ghz-rcs} where the superior accuracy of high-order discretizations is evident. Plots of the imaginary component of surface current for each order with mesh h0 are shown in Figures~\ref{subfig:almond-3ghz-surface-current-loworder} to \ref{subfig:almond-3ghz-surface-current-highorder} which illustrate the smoothness in the solution obtained at higher orders. 

\begin{figure}[htp]
	\centering
	\begin{subfigure}[b]{0.6\textwidth}
		\includegraphics[width=1.0\textwidth]{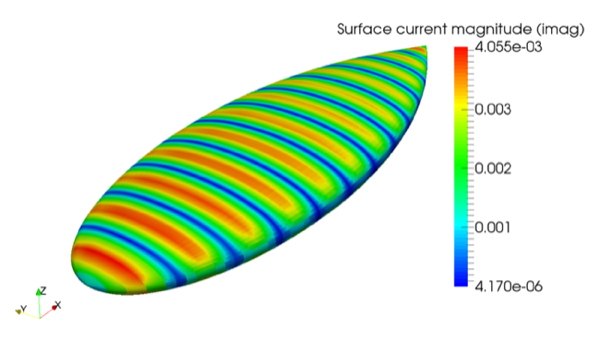}
		\caption{$(1,0)\nospacetimes(0,1)$.}
		\label{subfig:almond-7ghz-surface-current-loworder}
	\end{subfigure}
	
	\begin{subfigure}[b]{0.6\textwidth}
		\includegraphics[width=1.0\textwidth]{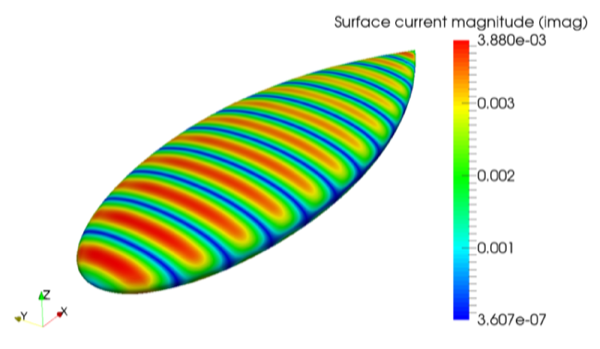}
		\caption{$(2,1)\nospacetimes(1,2)$.}
		\label{subfig:almond-7ghz-surface-current-midorder}  
	\end{subfigure}

	\begin{subfigure}[b]{0.6\textwidth}
		\includegraphics[width=1.0\textwidth]{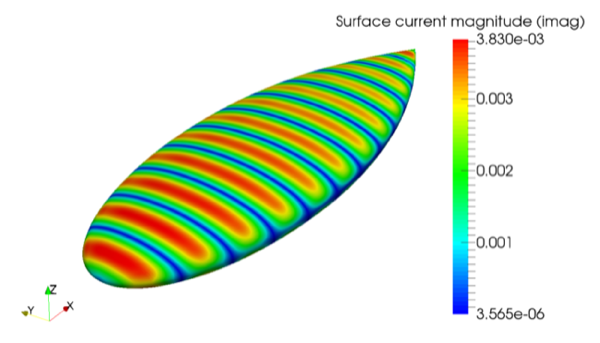}
		\caption{$(3,2)\nospacetimes(2,3)$.}
		\label{subfig:almond-7ghz-surface-current-highorder}  
	\end{subfigure}
	\caption{Magnitude of imaginary component of surface current over the NASA almond geometry: vertically polarised planewave of 7GHz travelling in the positive $x-$direction, mesh h2.}
	\label{fig:almond-7ghz-surface-current}
\end{figure}

Finally, we consider the 7GHz case where RCS plots for mesh h2 are illustrated in Figures~\ref{fig:almond-7ghz-hh-rcs} and   \ref{fig:almond-7ghz-vv-rcs} for HH and VV polarization respectively demonstrating good agreement with the numerical reference solution.  At this frequency large errors were encountered for meshes h0 and h1 necessitating the use of mesh h2. Plots of the imaginary component of surface current for each order with mesh h2 are shown in Figures~\ref{fig:almond-7ghz-surface-current}.

\subsection{Integrated CAD and electromagnetic scattering analysis}
\label{sec:direct_integration}

\begin{figure}[htp]
	\centering
	\begin{subfigure}[b]{0.45\textwidth}
		\includegraphics[width=1.0\textwidth]{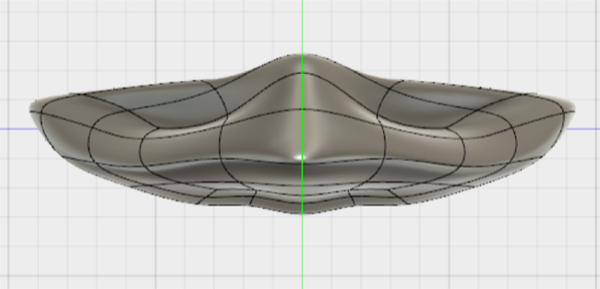}
		\caption{Front view ($y\operatorname{-}z$ plane).}
		\label{subfig:concept-stealth-front-view}
	\end{subfigure}
	\begin{subfigure}[b]{0.45\textwidth}
		\includegraphics[width=1.0\textwidth]{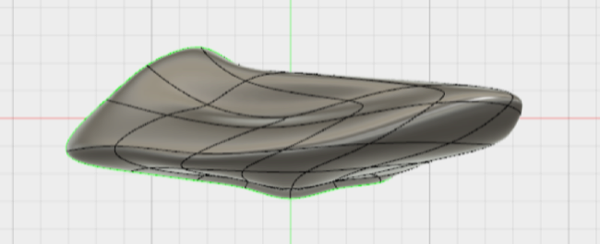}
		\caption{Side view ($x\operatorname{-}z$ plane).}
		\label{subfig:concept-stealth-sideview}
	\end{subfigure}

	\begin{subfigure}[b]{0.8\textwidth}
		\includegraphics[width=1.0\textwidth]{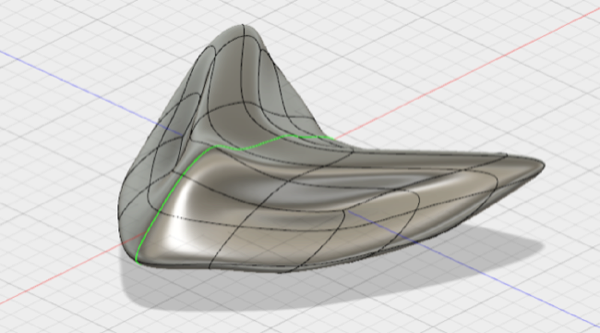}
		\caption{Perspective view.}
		\label{subfig:concept-stealth-perspective-view}  
	\end{subfigure}
	\caption{Concept model generated in Autodesk\textsuperscript{\textregistered} Fusion 360\texttrademark used for RCS analysis.}
	\label{fig:concept-stealth-cad}
\end{figure}

We now demonstrate the ability of our approach to perform electromagnetic scattering analysis directly on CAD generated models.  Figure~\ref{fig:concept-stealth-cad} illustrates a concept model generated in Autodesk\textsuperscript{\textregistered} Fusion 360\texttrademark~ which includes T-spline functionality capable of producing smooth, watertight surfaces.  The model is composed of six bicubic NURBS 
surfaces consisting of 1,178 control points and 384 elements.  By exporting this model as a STEP file which preserves all NURBS data structures and making use of the OpenCascade library, a compatible B-spline discretization is generated directly from this NURBS geometry model.  We envisage a scenario where our implementation could be included directly with a CAD software library thereby eliminating this STEP file export procedure.  The size of the bounding box for this model is given by $(\Delta x, \Delta y, \Delta z) = (82.3, 93.1, 27.5)$.

\begin{figure}[htp]
	\centering
	\begin{subfigure}[b]{0.8\textwidth}
		\includegraphics[width=1.0\textwidth]{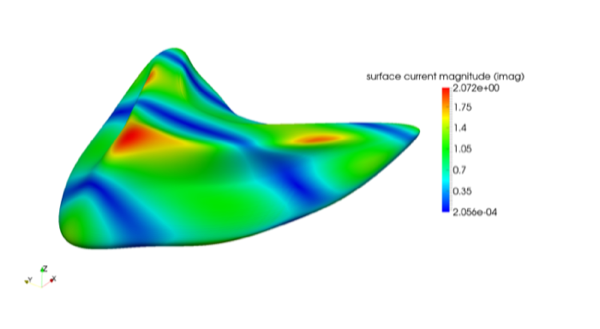}
		\caption{$|\mathbf{J}_{i}|$}
		\label{subfig:concept-stealth-surface-current}
	\end{subfigure}

	\begin{subfigure}[b]{0.8\textwidth}
		\includegraphics[width=1.0\textwidth]{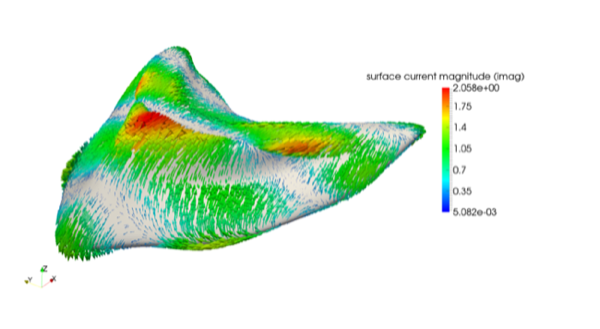}
		\caption{$\mathbf{J}_{i}$}
		\label{subfig:concept-stealth-surface-vector}  
	\end{subfigure}
	\caption{An example surface current profiles (imaginary) for the concept model shown in Figure~\ref{fig:concept-stealth-cad}. The plane wave is polarised in the $z$ direction and travelling in the positive $x$ direction with $ka = 9.31$.}
	\label{fig:concept-stealth-surface-current}
\end{figure}

\begin{figure}[htp]
  \centering
\includegraphics[width=0.8\textwidth]{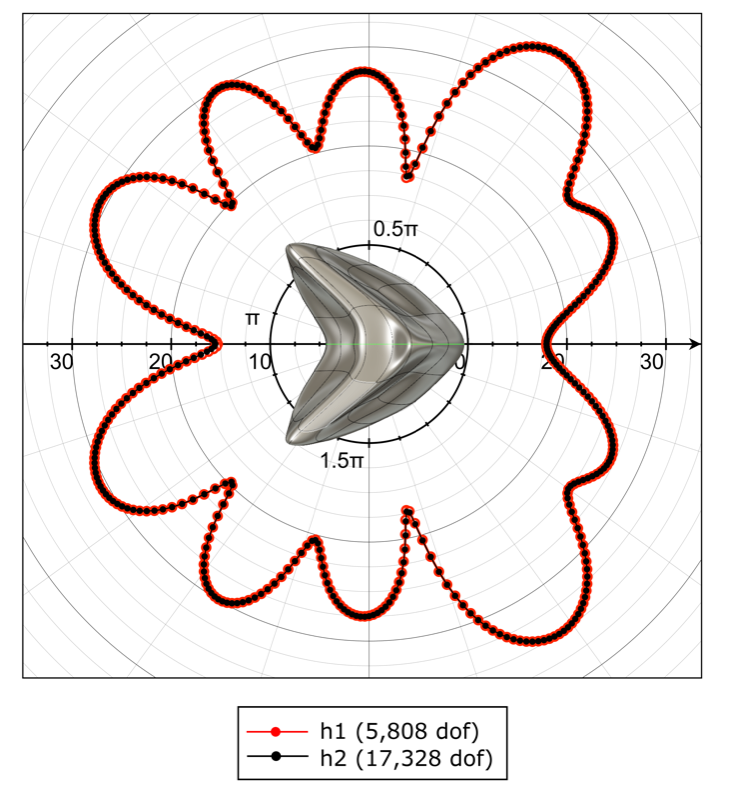}
  \caption{The computed radar cross section profile for the concept model illustrated in Figure~\ref{fig:concept-stealth-cad} with a normalised wavenumber $ka = 9.31$. }
  \label{fig:concept-stealth-rcs}  
\end{figure}

\begin{figure}[htp]
	\centering
	\begin{subfigure}[b]{0.8\textwidth}
		\includegraphics[width=1.0\textwidth]{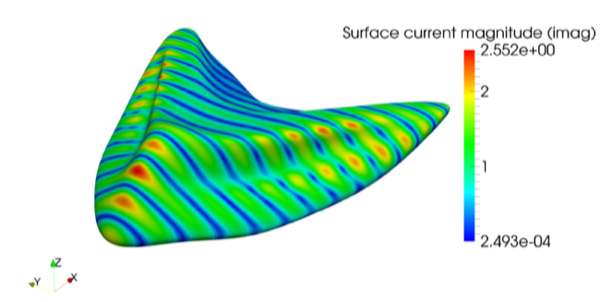}
		\caption{$|\mathbf{J}_{i}|$}
		\label{subfig:concept-stealth-surface-current-high-k}
	\end{subfigure}
	
	\begin{subfigure}[b]{0.8\textwidth}
		\includegraphics[width=1.0\textwidth]{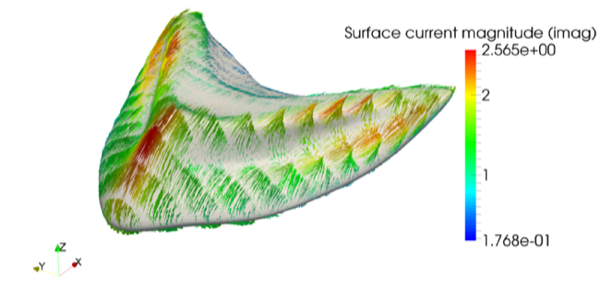}
		\caption{$\mathbf{J}_{i}$}
		\label{subfig:concept-stealth-surface-vector-high-k}  
	\end{subfigure}
	\caption{An example surface current profiles (imaginary) for the concept model shown in Figure~\ref{fig:concept-stealth-cad}. The plane wave is polarised in the $z$ direction and travelling in the positive $x$ direction with $ka = 46.55$.}
	\label{fig:concept-stealth-surface-current}
\end{figure}

\begin{figure}[htp]
	\centering
	\includegraphics[width=0.8\textwidth]{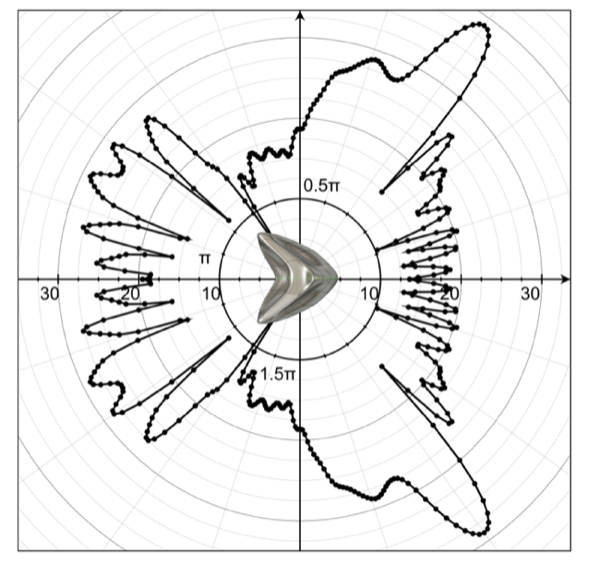}
	\caption{The computed radar cross section profile for the concept model illustrated in Figure~\ref{fig:concept-stealth-cad} with a normalised wavenumber $ka = 46.55$. }
	\label{fig:concept-stealth-rcs-high-k}  
\end{figure}

RCS values are computed over the $x\operatorname{-}y$ plane in which the wave is polarised in the $z$-direction.  We first apply a normalised wavenumber $ka = 9.31$ and apply two levels of h-refinement (denoted by $h1$ and $h2$ respectively)  using compatible B-splines of order $(3,2)\nospacetimes(2,3)$ with normal $C^{0}$ continuity across patches.  The discretizations $h1$ and $h2$ consist of 5,808 and 17,328 degrees of freedom respectively.  Plots of the imaginary component of surface current for $h2$ are shown in Figures~\ref{subfig:concept-stealth-surface-current} and \ref{subfig:concept-stealth-surface-vector} and RCS values are plotted in Figure~\ref{fig:concept-stealth-rcs}.  We also compute RCS values for a higher normalised wavenumber of $ka = 46.55$ in which three levels of h-refinement are applied generating a discretisation with 58,800 degrees of freedom.  Surface current plots for this wavenumber are shown in Figures~\ref{subfig:concept-stealth-surface-current-high-k} and \ref{subfig:concept-stealth-surface-vector-high-k} and RCS values are plotted in Figure~\ref{fig:concept-stealth-rcs-high-k}.


We use this example to demonstrate how our approach exhibits a tight link between computational design and analysis by using a common data model that provides the necessary geometry and analysis discretizations.  The requirement for surface meshing is bypassed and the use of high order B-spline discretizations provides superior accuracy over conventional discretization approaches. 

\section{Conclusion}
\label{sec:conclusion}

We have outlined an isogeometric boundary element method (method of moments) that utilises a common model to discretise both the geometry and analysis fields for electromagnetic scattering analysis.  Our approach uses Non-Uniform Rational B-Splines (NURBS) to represent the surface geometry and compatible B-splines as basis for electromagnetic analysis.  We have detailed the construction of compatible B-splines from a given NURBS discretization  that provide a div-conforming or curl-conforming surface vector basis and described how such spline-based discretizations can be used as a basis for the  electric/magnetic field integral equations.  We verified our approach through the Mie series solution that provides a closed-form solution for electromagnetic scattering over a perfectly electrically conducting sphere and utilised experimental and numerical reference data for the well-known NASA almond geometry to verify radar cross section calculations.  Finally, we demonstrated how our approach can be used to perform electromagnetic scattering analysis directly on geometry models generated using modern CAD software showcasing the ability of our approach to fully integrate CAD and analysis technologies.

\section*{Acknowledgements}
\label{sec:aknowledgements}
J.A. Evans was partially supported by the Air Force Office of Scientific Research under Grant No. FA9550-14-1-0113.

\begin{appendices}

\section{MFIE: compatible B-spline discretization}
\label{app:mfie}

In a similar manner to the electric field integral equation, the magnetic field integral equation is first derived by substituting the expression for the total magnetic field given by 
\begin{equation}
  \mathbf{H} = \mathbf{H}^{i} + \mathbf{H}^s.\label{eq:Htotal}
\end{equation}
into the PEC condition of 
\begin{equation}
\mathbf{n} \times \mathbf{H} = \mathbf{J} \label{eq:PECcondition2} 
\end{equation}
to arrive at
\begin{equation}
  \label{eq:mfie_1}
  \mathbf{n} \times \mathbf{H}^i  = \mathbf{J} - \mathbf{n} \times \mathbf{H}^s
\end{equation}
with the scattered magnetic field given by the quantity 
\begin{equation}
  \label{eq:mfie2}
   \mathbf{H}^s = \nabla \times \mathbf{A}
\end{equation}
allowing (\ref{eq:mfie_1}) to be rewritten as
\begin{equation}
  \label{eq:mfie3}
   \mathbf{n} \times \mathbf{H}^i  = \mathbf{J} - \mathbf{n} \times \int_\Gamma \nabla \times \mathbf{J} \,\frac{e^{-jkr}}{4 \pi r} \,\mathrm{d}\Gamma.
\end{equation}
Defining the linear operator
\begin{equation}
  \label{eq:mfie_linearoperator}
  L^H(\mathbf{u}) =  \mathbf{u} - \mathbf{n} \times \int_\Gamma \nabla \times \mathbf{u} \,\frac{e^{-jkr}}{4 \pi r} \,\mathrm{d}\Gamma
\end{equation}
and a forcing function $\mathbf{g} = \mathbf{n} \times \mathbf{H}^i$, we write the Galerkin formulation of the magnetic field integral equation as:

given $\mathbf{g}$,  find $\mathbf{J} \in H^{-\frac{1}{2}}(\mathrm{curl}_\Gamma, \Gamma)$ such that 
\begin{equation}
  \label{eq:mfie_forcingterm}
  \langle \mathbf{w}, L^H(\mathbf{J}) \rangle = \langle \mathbf{w}, \mathbf{g}  \rangle \quad \forall \mathbf{w} \in  H^{-\frac{1}{2}}(\mathrm{\mathrm{curl}}_\Gamma, \Gamma).
\end{equation}
Defining finite dimensional subspaces $\mathbf{w}_h, \mathbf{J}_h \in H^{-\frac{1}{2}}(\mathrm{curl}_\Gamma, \Gamma)$  as
\begin{align}
\mathbf{w}_h &= \sum_A^{N_{b}} \mathbf{N}^{\mathrm{curl}}_A w_A \label{eq:weight_discretization_curl}\\
 \mathbf{J}_h &=  \sum_A^{N_{b}} \mathbf{N}^{\mathrm{curl}}_A j_A \label{eq:current_discretization_curl}
\end{align}
where $\{\mathbf{N}^{\mathrm{curl}}_A\}_{A=1}^{N_b}$ is a set of curl-conforming surface vector B-spline basis functions , the system of equations for the magnetic field integral equation can be written as
\begin{equation}
  \label{eq:system_equations_mfie}
  \mathbf{Y}_{AB} \mathbf{J}_B = \mathbf{g}_A
\end{equation}
where,  by employing the identity $\nabla \times (\phi \mathbf{v}) = \nabla \phi \times \mathbf{v} + \phi \nabla \times \mathbf{v}$, applying a limiting process to the integral and noting that $\mathbf{N}_A^{\mathrm{div}} = -\mathbf{n} \times \mathbf{N}_A^{\mathrm{curl}}$, 
\begin{equation}
  \label{eq:mfie_matrix}
  \mathbf{Y}_{AB} = \frac{1}{2} \int_{\Gamma_x} \mathbf{N}_A^{\mathrm{curl}} \cdot \mathbf{N}_B^{\mathrm{curl}}\,\mathrm{d}\Gamma +  \int_{\Gamma_x} \mathbf{N}_A^{\mathrm{div}} \cdot \left(  \int_{\Gamma_y} \nabla G \times \mathbf{N}_A^{\mathrm{curl}}\,\mathrm{d}\Gamma  \right) \,\mathrm{d}\Gamma
\end{equation}
where
\begin{equation}
  \label{eq:gradient_greensfunction}
  \nabla G = -\frac{e^{-jkr}}{4 \pi r} \left( \frac{1}{r} + jk \right) \mathbf{r}
\end{equation}
with $\mathbf{r} := \mathbf{y} - \mathbf{x}$ and the factor of $1/2$ arises from the limiting process. Similarly, the forcing vector components are given by 
\begin{align}
  \label{eq:mfie_forcing_disc}
 \mathbf{g}_A &= \int_{\Gamma_x} \mathbf{N}_A^{\mathrm{curl}} \cdot (\mathbf{n} \times \mathbf{H}^i)  \,\mathrm{d}\Gamma\\
&=\int_{\Gamma_x} \mathbf{N}_A^{\mathrm{div}} \cdot \mathbf{H}^i  \,\mathrm{d}\Gamma.
\end{align}
As before, the vector $\mathbf{J}_B$ represents a vector of unknown  surface current density coefficients.

%

\section{PEC sphere - additional results}
\label{app:sec:pec-additional-results}

\begin{figure}[H]
    \centering
    \begin{subfigure}[b]{0.6\textwidth}
      \centering
    \includegraphics[width=1.0\textwidth]{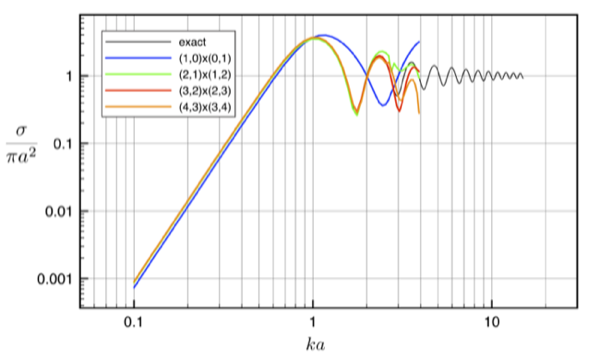}
      \caption{h0}  
      \label{fig:rcs-sphere-h0}
    \end{subfigure}

    \begin{subfigure}[b]{0.6\textwidth}
      \centering
    \includegraphics[width=1.0\textwidth]{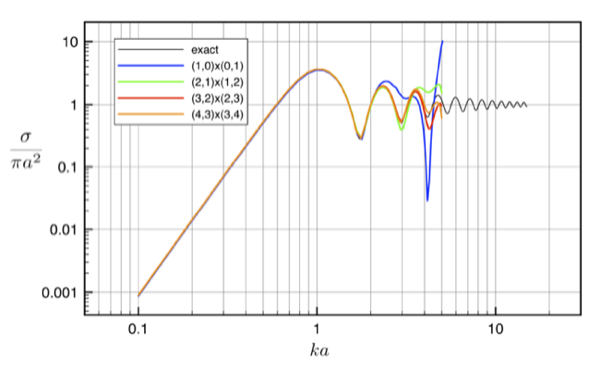}
      \caption{h1}  
      \label{fig:rcs-sphere-h1}
    \end{subfigure}

    \begin{subfigure}[b]{0.6\textwidth}
      \centering
    \includegraphics[width=1.0\textwidth]{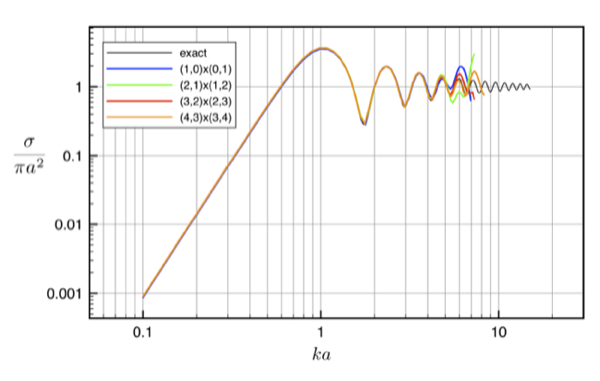}
      \caption{h2}  
      \label{fig:rcs-sphere-h2}
    \end{subfigure}
    
    \caption{Normalised RCS values for a PEC sphere computed for increasing wavenumber with div-conforming B-splines of varying degree.}
    \label{fig:sphere-rcs-h0-to-h2}
\end{figure}

%

\section{NASA almond geometry parameterization}
\label{app:sec:nasa-geometry}

Denoting the length of the almond geometry as $L=0.2524m$, the surface of the NASA almond geometry is defined in terms of parametric coordinates $(s,t)$ as 
\begin{align}
\begin{pmatrix}
x\\
y\\
z
\end{pmatrix}
&=
\begin{pmatrix}
Lt\\
0.193333L \sqrt{1 - \left( \frac{t}{0.416667}\right)^2 } \cos s\\
0.064444L \sqrt{1 - \left( \frac{t}{0.416667}\right)^2 	} \sin s
\end{pmatrix}\\
&\textrm{for} -\pi< s< \pi, -0.41667 < t < 0\notag
\end{align} 
and
\begin{align}
\begin{pmatrix}
x\\
y\\
z
\end{pmatrix}
&=
\begin{pmatrix}
Lt\\
4.83345L \left[\sqrt{1 - \left( \frac{t}{2.08335}\right)^2}-0.96 \right] \cos s\\
1.61115L \left[\sqrt{1 - \left( \frac{t}{2.08335}\right)^2} - 0.96\right] \sin s
\end{pmatrix}\\
&\textrm{for} -\pi< s< \pi, 0 < t < 0.58333.\notag
\end{align} 

\end{appendices}

\bibliographystyle{elsarticle-num} 
\bibliography{bibliography}

\end{document}